\documentclass[11pt]{article}
\usepackage{graphicx}
\usepackage{amscd}
\usepackage{amsmath}
\usepackage{amsfonts}
\usepackage{amssymb}

\newcommand{\be}{\begin{eqnarray*}}
\newcommand{\ee}{\end{eqnarray*}}

\newcommand{\al}{\alpha}

\newcommand{\f}[1]{{\frak #1}}

\newcommand{\lon}{\longrightarrow}
\newcommand{\complex}{{\Bbb C}}

\newcommand{\half}{\frac{1}{2}}
\newcommand{\third}{\frac{1}{3}}

\newcommand{\cinf}{C^{\infty}}

\newtheorem{thm}{Theorem}[section]
\newtheorem{lem}[thm]{Lemma}
\newtheorem{cor}[thm]{Corollary}
\newtheorem{pro}[thm]{Proposition}

\newtheorem{rmk}[thm]{Remark}
\newcommand{\pf}{\noindent{\bf Proof.}\ }
\newcommand{\qed}{\begin{flushright} $\Box$\ \ \ \ \  \end{flushright}\\\\}
\newtheorem{defi}[thm]{Definition}

\newcommand{\frakg}{{\frak g}}
\newcommand{\frakh}{{\frak h}}
\newcommand{\frakm}{{\frak m}}

\newcommand{\calx}{{\cal X}}
\newcommand{\call}{{\cal L}}
\newcommand{\calo}{{\cal O}}
\newcommand{\gm}{\Gamma }

\newcommand{\ot}{\otimes}
\newcommand{\Alt}{\mbox{Alt}}
\newcommand{\smalcirc}{\mbox{\tiny{$\circ $}}}
\newcommand{\h}{{\frak h}}
\newcommand{\g}{{\frak g}}
\newcommand{\xa}{\xi_{\al}}
\newcommand{\ax}{\xi_{-\al}}
\newcommand{\calp}{{\cal P}}
\newcommand{\calu}{{\cal U}}

\newcommand{\cala}{{\cal A}}
\newcommand{\del}{\partial}

\newcommand{\et}{\epsilon \downarrow 0}
\newcommand{\ct}{c\to }
\newcommand{\CM}{\mbox{Calogero-Moser }}
\newcommand{\cm}{\mbox{Calogero-Moser }}
\newcommand{\hp}{\frakh^{\perp}}
\newcommand{\ppi}{\pi}
\newcommand{\HH}{{\cal H}}
\newcommand{\DD}{H}
\newcommand{\ra}{\rangle}
\newcommand{\la}{\langle}
\newcommand{\xx}{\tilde{x}}

\def\qed{\hfill ~\vrule height6pt width6pt depth0pt\bigskip}

\def\fg{{\frak{g}}}
\def\fh{{\frak{h}}}
\def\dx{\dot{x}}
\def\wl{\widetilde{L}}
\def\wr{\widetilde{R}}

\def\wphi{\widetilde{\phi}}
\def\wpsi{\widetilde{\psi}}

\begin{document}

\title{{\bf Integrable spin Calogero-Moser systems}}

\author{Luen-Chau Li\\
{\sf email: luenli@math.psu.edu }\\
AND\\
Ping Xu
         \thanks{ Research partially supported by NSF
       grant DMS00-72171. }\\
{\sf email: ping@math.psu.edu }\\
        Department of Mathematics\\
         Pennsylvania State University \\
         University Park, PA 16802, USA}

\date{}

\maketitle

\begin{abstract}
We introduce spin Calogero-Moser systems associated with root systems
of simple Lie algebras and give the associated Lax representations
 (with spectral parameter) and fundamental Poisson
 bracket relations. The associated integrable models (called
integrable spin Calogero-Moser systems in the paper) and their
Lax pairs are then obtained via Poisson reduction
and gauge transformations. For  Lie algebras of $A_{n}$-type,
this new class of  integrable systems includes the usual Calogero-Moser systems 
as subsystems. Our  method is guided by a general framework which
we develop here using dynamical Lie algebroids.
\end{abstract}

\section{Introduction}

Calogero-Moser type systems are Hamiltonian systems with very rich structures.
After the pioneering work of Calogero and Moser \cite{C} \cite{M}, 
many generalizations have been proposed.  Olshanetsky and Perelomov
\cite{OP}, for example,
introduced Calogero-Moser models associated with root systems of simple Lie
algebras (for recent work, see for example, \cite{BCS} and \cite{DP}).
  On the other hand, a rational ${\frak{sl}}(N, \complex) $
 spin Calogero-Moser system
was introduced by Gibbons and Hermsen \cite{GH}.  As in the spinless case,
trigonometric and elliptic versions of this generalization also exist.  In
recent years, these models and their variants have received considerable
attention due to their relevance in a number of areas.  In the original work
of Gibbons and Hermsen, and in the paper \cite{KBBT}, for example, 
the ${\frak{sl}}(N, \complex)$ spin 
systems considered by the authors are related to certain special solutions of
integrable partial differential equations.  In a completely different area,
an approach to study the joint distribution of energy eigenvalues of a
Hamiltonian was initiated by Pechukas \cite{P} and continued by Yukawa 
\cite{Y} and a number of other researchers (see, for example,
\cite{HM} and the references 
therein).  In this so-called level dynamics approach in random matrix
theory,  spin Calogero-Moser systems appear naturally.  As a matter
of fact, they provide the starting point of the ensuing analysis.  For a
recent connection between $SU(2)$ Yang-Mills mechanics and a version of
the rational model embedded in an external field, we refer the reader to
\cite{KM}.  At this juncture, we should perhaps warn the reader over possible
confusion with the term spin Calogero-Moser systems.  Indeed, there are
many different versions of this kind of generalization and yet the same term 
is used to describe these different systems.  For example, in \cite{GH} and
\cite{KBBT},
the authors were actually restricting themselves to a special symplectic leaf
of an  underlying Poisson manifold.  On the other hand, in \cite{P} and 
\cite{Y}, the spin variables are in the space of skew-Hermitian matrices.
  In this regard,
the reader can consult \cite{W} and \cite{Po}
 for further examples in addition to those
mentioned above.  See also Remark \ref{rmk:4.2} (2)  in Section 4.

     In \cite{BAB1}, the authors considered the rational
 ${\frak{sl}}(N , \complex )$
 spin Calogero-Moser system.  Without restricting themselves to a special symplectic leaf
as in \cite{KBBT}, they obtained a St. Petersburg type formula
 for the ${\frak{sl}}(N , \complex )$
model, i.e., the so-called fundamental Poisson bracket relation (FPR) between 
the elements of an associated Lax operator $L(z)$.  However, what they
found was rather unusual.  First of all, there are the usual kind of terms
in the FPR, but now an r-matrix depending on phase space variables is
involved.  Then there is an anomalous term whose presence is an obstruction
to integrability.  By this, we mean that the quantities $\mbox{tr}(L( z)^{n})$
do not Poisson commute unless we restrict to the submanifold $\Sigma$ where
the anomalous term vanishes.  If $\Sigma$ were a Poisson submanifold of the
underlying Poisson manifold, the corresponding subsystem would have a 
natural collection of Poisson commuting integrals, but unfortunately this
is not the case.  For the trigonometric and elliptic 
${\frak{sl}}(N , \complex )$ systems, similar
formulas were obtained in \cite{BAB2}.

     Our present work has its origin in an attempt to understand conceptually
the group theoretic/geometric meaning of the wonderful but mysterious
calculations in \cite{BAB1} and \cite{BAB2}.
  As was pointed out in a later paper by the same authors \cite{ABB},
the r-matrices which appear in their earlier work do satisfy a closed-
form equation, the so-called classical dynamical Yang-Baxter equation
(CDYBE) \cite{F}.  CDYBE is an important differential-functional equation
introduced by Felder in his work on conformal field theory \cite{F}.
For simple Lie algebras and Kac-Moody algebras, the classification of solutions of
this equation (under certain conditions) was obtained by Etingof and
Varchenko \cite{EV}.  On the other hand, dynamical r-matrices are intimately
related to coboundary dynamical Poisson groupoids \cite{EV} and coboundary Lie
bialgebroids \cite{BK-S}.
  This relation is analogous to the more familiar one which
exists between constant r-matrices, Poisson groups and Lie bialgebras \cite{D}.
Consequently, it is plausible that the calculations in \cite{BAB1} and
\cite{BAB2} are
connected with Lie algebroids, and as it turns out, this is indeed the case. In
this connection, let us recall that in integrable systems theory, one of the 
powerful means to show that a Hamiltonian system is integrable (in some sense)
is to realize the system in the r-matrix scheme for constant r-matrices 
(see \cite{FT} \cite{RSTS} and the references therein). 
 For the ${\frak{sl}}(N , \complex )$ spin Calogero-Moser systems,
 we have found an analog of the realization picture, using 
Lie algebroids associated with dynamical r-matrices.  Indeed, along the way,
it became clear that one can introduce spin Calogero-Moser systems associated
with root systems of simple Lie algebras.  These spin systems are naturally
associated with the dynamical r-matrices with spectral parameter in \cite{EV}.
Furthermore, there is a unified way to construct the realization maps for such
systems.  However, as in the ${\frak{sl}}(N , \complex )$ case, 
there is an obstruction for the 
natural functions to Poisson commute.  Nevertheless, the underlying structures
of the spin systems permits the construction of associated integrable models,
via Poisson reduction \cite{MR} and the idea of gauge 
transformations \cite{BAB1}.
More precisely, the Hamiltonians of the spin Calogero-Moser systems are 
invariant
 under a natural canonical action of a Cartan subgroup of the underlying
simple Lie group.  In addition, the obstruction to integrability vanishes on a
fiber of the equivariant momentum map.  Hence we can apply Poisson reduction
to obtain the integrable models on reduced Poisson manifolds.  
We shall call the systems in this {\em new} class of integrable
models {\em integrable spin Calogero-Moser systems}.

     We now describe the contents of the paper.  In Section 2, we assemble
a number of basic facts which will be used in the paper.  In Section 3, we
consider realization of Hamiltonian systems in dynamical Lie algebroids.  
The reader should note the substantial difference between our present case
and the more familiar case of realization in Lie algebras equipped with an
R-bracket.  The difference lies in the fact that in the present case, the
natural functions to consider are functions which do not Poisson commute
on the dual of the dynamical Lie algebroid.  Consequently, we do not get
integrable systems to start with.  However, when the realization map is
an equivariant map, then under suitable assumptions, we show that Poisson
reduction can be used to produce integrable flows with natural family of
conserved quantities in involution on reduced phase spaces.  An important
idea which we employ in this connection is that of a gauge transformation
of a Lax operator which we learned from the paper \cite{BAB1}.  As the reader
will see, this device not only allows us to write down the equation of motion
in Lax pair form on the reduced Poisson manifold.  It also enables us to
establish involution of the induced functions.  In Section 4, we introduce
spin Calogero-Moser systems associated with root systems of simple Lie
algebras.  Our first step in this section is the construction of dynamical
Lie algebroids, starting from the classical dynamical r-matrices with
spectral parameter.  Then we show that the putative Poisson manifold
underlying the spin Calogero-Moser systems admit realizations in the dynamical
Lie algebroids constructed earlier.  Here, the realization maps are natural
in the sense that the corresponding dual maps are morphisms of Lie algebroids.
In \cite{EV}, Etingof and Varchenko obtained a classification of classical
dynamical r-matrices with spectral parameter for simple Lie algebras.  They
obtained canonical forms of the three types of dynamical r-matrices (rational,
trigonometric, and elliptic).  For each of these canonical forms, we can
use the corresponding realization map to construct the associated spin
Calogero-Moser systems.  In Section 5, we carry out the reduction procedure
to the spin systems to obtain the associated integrable models.  Here, the
main task is to construct an equivariant map from an open dense subset of
the Poisson manifold of the spin systems to the Cartan subgroup.  Using
this map, we can define gauge transformations of the Lax operators for
the spin systems.  These gauge transforms are invariant under the natural
action of the Cartan subgroup and hence descend to the reduced Poisson
manifold.  In this way, we can obtain the Lax equations for the reduced
Hamiltonian systems and establish an involution theorem for the induced
functions. 
Furthermore, this gives rise to spectral curves
 which are preserved by the Hamiltonian flows.
     Some of the results in Section 4 of the present work have been announced
in \cite{LX}.


{\bf Acknowledgments.}  
We  would like to thank several institutions
for their hospitality while work on this project was being done:
MSRI (Li),  and Max-Planck Institut (Li and Xu). The first author
thanks the organizers, Pavel Bleher and Alexander Its, of the special
semester in Random matrix models and their
Applications held at MSRI in Spring 1999 for
hospitality during his stay there.  We  also  wish to thank
Jean Avan, Pavel Etingof, Eyal  Markman and Serge Parmentier 
 for discussions.

\section{Dynamical Lie algebroids}

\label{sec-algebroids}

In this section, we recall some basic facts.
Most of the  material is standard, which is presented 
here for the reader's convenience.

A Lie algebroid over a manifold $M$ may be thought of as a
``generalized tangent bundle'' to $M$.  Here is the definition (see
\cite{ma:lie} \cite{CW} for more details on the theory).

\begin{defi}
\label{dfn-algebroid}
A Lie algebroid over a manifold $M$ is a vector bundle $A $ over
$M$ equipped with a Lie
algebra structure $[\cdot , \cdot ]$ on its space of sections and a bundle map
$a : A \rightarrow TM$ (called the {\em anchor}) such that
\begin{enumerate}
\item the bundle map $a: A\lon TM$
 induces a Lie algebra homomorphism (also denoted by
 $a $) from sections of $A $ to vector fields on $M$;
\item for any $X, Y\in \gm (A)$ and $f\in C^{\infty}(M)$, the
identity
$$ [X , fY ] = f[X ,Y ]+ (a(X)f)Y$$ 
holds.
\end{enumerate}
\end{defi}

Examples of Lie algebroids include the usual Lie algebras,
Lie algebra bundles,  tangent bundles of smooth manifolds,
 and integrable distributions  on smooth manifolds.
If $A$ is finite-dimensional, the standard local coordinates
 on $A $ are of  the form $(q,\lambda)$,
where the $q_{i}$'s are coordinates on the base $M$ and the $\lambda
_{i}$'s are linear coordinates on the fibers, associated with a basis
$X_{i}$ of sections of the Lie algebroid.  In terms of such
coordinates, the bracket and anchor are given by: 
\begin{equation}
\label{eq:str-constant}
[X_{i}, X_{j}]= \sum c_{ij}^{k} X_{k}, \ \ \mbox{ and }
a (X_{i})=\sum a_{ij}\frac{\del }{\del q_{j}},
\end{equation}
 where  $c_{ij}^{k}$ and $a_{ij}$ are ``structure
functions'' lying in $\cinf(M)$.

The dual bundle $A^{*}$   of $A$ carries a natural
Poisson structure, called the Lie-Poisson structure
\cite{CDW}.  To describe this structure, it suffices to give
the Poisson brackets of a class of functions whose differentials span
the cotangent space at each point of $A ^{*}$.  Such a class is
given by functions which are affine on fibres.  The functions which
are constant on fibres are just the functions on $M$, lifted to $A ^{*}$ 
via the bundle projection. On the other hand,
functions which are  linear on fibres
may be identified with the sections of $A $.
This is because   for any $X\in \gm (A)$, we can define
$l_{X}\in C^{\infty} (A^* )$  by $l_{X} (\xi )=
<\xi , X>,\ \ \forall \xi\in A^*$.  If $f$ and $g$ are
functions on $M$, and $X $ and $Y $ are sections of $A $,
the Lie-Poisson structure is characterized by the following
 bracket relations:
 $$ \{f,g\}=0, ~\{f,  l_{X}\}= a(X) ( f ),
\mbox{ and }\{l_{X} , l_{X} \}= l_{[X, Y ]} \enspace .$$
For the finite-dimensional case, corresponding to
 standard coordinates $(q,\lambda )$ on $A $, we may
introduce dual coordinates $(q,\mu )$ on $\cala ^{*}$.  In terms of
such coordinates and the structure functions introduced 
 in Equation (\ref{eq:str-constant}),
 the Poisson bracket relations on $A ^{*}$ are 
$$\{q_{i},q_{j}\}=0, \ \  \{\mu _{i},\mu _{j}\}=c_{ij}^{k}\mu _{k},
\ \ \mbox{ and } \{q_{i},\mu _{j}\}=a_{ji}.$$

The Poisson structure on  $A^*$ generalizes
the usual Lie-Poisson structure on the dual of a
Lie algebra. Namely, if $A $ 
is a Lie algebra $\frakg$, the Poisson structure on its dual is
the standard Lie-Poisson structure on $\frakg^*$. On the other hand,
when $A =TM$ is  equipped with the standard Lie algebroid structure,
 the Poisson structure on its dual is just the usual
cotangent bundle symplectic  structure.
 Another interesting example, which we need in this paper,
 is the following\\\\
{\bf Example 2.1}
Let   $A=TM\times \frakg$  be
 equipped with the standard  product
Lie algebroid structure; namely, the anchor is the projection
map onto the first factor and the bracket on sections
is given by
\begin{equation}
\label{gb}
  [(X , \xi), \  (Y ,\eta )] =( [ X,Y] , \  [ \xi, \eta] + L_X \eta -L_Y \xi),
   \, \,\, X ,Y \in  \gm (TM), \, \,  \xi, \eta \in C^{\infty}(M, \f{g}),
\end{equation}
where the bracket of two vector fields is the usual bracket and the bracket
$[ \xi, \eta]$ is the pointwise bracket.
 Then clearly,  $A^*$ is the Poisson manifold direct product
 $T^* M\times \frakg^*$. In other words, the bracket between
functions on $T^* M$ is the canonical one on $T^* M$,
the bracket between
functions on $\frakg^*$ is the Lie-Poisson bracket, and the mixed
term  bracket between functions on $T^* M$ and $\frakg^*$ is zero.\\\\



In the rest of the section,
let $\frakg$ be  a Lie algebra  and $\frakh$ an Abelian
Lie subalgebra of $\frakg$.
 Consider   $T^*\frakh^{*} \times \frakg^*$ as
a vector bundle over $\frakh^{*}$, and define  a bundle map
$a_{*}: T^*\frakh^{*} \times \frakg^*\lon T\frakh^{*}$
by 
\begin{equation}
a_{*}(q, p, \xi )=(q, i^{*} \xi ), \ \  \forall q\in \frakh^{*} ,
p\in \frakh , \ \mbox{ and }\  \xi\in \frakg^*,
\end{equation}
 where $i: \frakh\lon \frakg$ is the inclusion map.
If  $R$ is   a map from $\frakh^{*}$ to
 $\call (\frakg^* , \frakg )$ (the space of linear maps from
$\frakg^*$ to $\frakg$),
we define a bracket on $\gm (T^*\frakh^{*} \times \frakg^* )$ as follows.
For $\xi , \eta \in \frakg^*$  considered as constant
sections, $h\in \frakh $ considered as a constant
one form on $\frakh^{*}$, and $\omega , \theta \in \Omega^{1}(\frakh^{*} )$,
 define
\be
&& [\omega, \theta ]=0, \\
&&[h , \xi ]=ad_{h}^*\xi , \\ 
&&[\xi ,\eta ]=d\la R\xi , \eta \ra -ad^*_{R\xi }\eta +ad^*_{R\eta }\xi,
\ee
where $ad^*$ denotes the dual of $ad$: $\la ad^*_{X}\xi , Y\ra =
\la \xi, [X, Y]\ra , 
\forall X, Y\in  \frakg$ and $\xi \in \frakg^*$.
Then  extend  this  to a bracket $[\cdot , \cdot ]$
for all sections in  $\gm (T^*\frakh^{*} \times \frakg^* )$ by
the usual anchor condition.

The following proposition  can be  verified by a direct calculation.

\begin{pro}
\label{pro:R}
$(T^*\frakh^{*} \times \frakg^* ,  \ [\cdot , \cdot ])$ is a Lie algebroid
  with anchor map $a_{*} $ iff
\begin{enumerate}
\item The operator  $R$ is a  map from $\frakh^{*} $
to  $\call (\frakg^* , \frakg )^{\frakh}$,
 the space of $\frakh$-equivariant linear map from
$\frakg^*$ to $\frakg$  ($\frakh$ acts on $\frakg$ by
adjoint action and on $\frakg^*$ by coadjoint action);
\item   $R$ satisfies $R_{q}^* =-R_{q}$ for
each point $q\in \frakh^{*}$ (here, as well as in the sequel,
we denote by $R_{q}$   the linear map in $\call (\frakg^* , \frakg )$
obtained by evaluating $R$ at the point $q$);
\item For any $q \in \frakh^{*}$,
the linear map from $\frakg^* \otimes \frakg^* \lon \frakg$ 
defined by
\begin{equation}
\xi \otimes \eta  \lon
[R_q \xi , R_q \eta ]+R_q (ad^*_{R_q \xi }\eta -ad^*_{R_q \eta }\xi)
+ \calx_{i^* \xi} (q) (R \eta ) -\calx_{i^* \eta}(q) (R \xi ) +d\la R \xi , \eta \ra  (q)
\end{equation}
is independent of $q \in \frakh^{*}$, and is  $\frakg$-equivariant,
where $\frakg$ acts on $\frakg^* \otimes \frakg^*$ by coadjoint action
and on $\frakg$ by adjoint action. Here, as well as in the sequel, 
 $\calx_v$ for $v\in \frakh^*$ denotes 
the  operation of taking the  derivative with respect to
$q$  along the constant vector field defined by $v$.
\end{enumerate}
\end{pro}
Such a  Lie algebroid  $(T^*\frakh^{*} \times \frakg^* ,  \ [\cdot , \cdot ])$
will be called a {\em dynamical Lie algebroid}, and we shall
use this terminology throughout the paper.

\begin{rmk}
\label{rmk:2.1}
{\em If $\frakg$ is finite-dimensional, and
 $R_{q}=r(q)^{\#}$ (i.e., $\la R_{q}\xi, \eta 
\ra =\la r(q) , \ \xi  \ot \eta \ra, \ \xi, \eta\in \frakg^*$)
 for a map $r: \frakh^{*}  \lon \wedge^{2}\frakg $, it can be shown 
that $R$ satisfies  the conditions in Proposition \ref{pro:R}, iff
$r$ satisfies:
\begin{enumerate}
\item $r$ is $\frakh$-invariant, i.e.,
$[1\ot h +h\ot 1 , \ r(q )]=0, \ \forall q\in \frakh^*, \  h\in \frakh$;
\item  $\sum_{i} h_{i} \wedge \frac{\partial r}{ \partial q_{i}} 
+\half [r, r]$ is a constant
$(\wedge^{3} \frakg)^{\frakg}$-valued function  over $\frakh^{*}$,
where $[\cdot , \cdot ]$ is the Schouten bracket
on $\oplus \wedge^* \frakg$,
 $\{h_{1}, \cdots , h_{N}\}$  is a basis in  $\frakh$,  and
$(q_{1}, \cdots , q_{N})$ its  induced coordinate system  on $\frakh^*$. 
\end{enumerate}
In other words, $r$ is  a dynamical $r$-matrix in the sense of \cite{F} 
\cite{EV}.
 Indeed, $(T\frakh^{*} \times \frakg , \ T^*\frakh^{*} \times \frakg^* )$
is   a Lie bialgebroid \cite{BK-S}.
}
\end{rmk}

Next,  we assume that  $\frakg$ admits a non-degenerate ad-invariant pairing
$(\cdot , \cdot )$.  
If $I: \frakg^* \lon \frakg$
is the  induced  isomorphism,  then a straightforward
calculation yields that
\begin{equation}
I  (ad^*_{X}\xi ) =-[X, I\xi ], \ \ \forall X\in \frakg, \ \xi\in \frakg^*.
\end{equation} 
 Thus we have the following

\begin{cor}
\label{cor:R}
The operator $R:  \frakh^{*} \lon \call (\frakg^* , \frakg )$ defines
a Lie algebroid structure on  $T^*\frakh^{*} \times \frakg^*$
if  the condition (1)-(2) in Proposition \ref{pro:R} are satisfied,
 and if   $R$ satisfies the modified dynamical Yang-Baxter equation (mDYBE):
\begin{equation}
\label{eq:mdybe}
[R\xi , R\eta ]-  R(I^{-1}[R\xi , I\eta ]+I^{-1}[I\xi , R\eta ] )
+ \calx_{i^* \xi}(R\eta ) -\calx_{i^* \eta}(R\xi ) +d\la R\xi , \eta\ra  
=c[I\xi , I\eta ],  \ \ \forall \xi, \eta \in \frakg^* ,
\end{equation}
for some constant $c$.
\end{cor}

\section{Realization of Hamiltonian systems in dynamical Lie
algebroids}

Throughout this section, let $T^*\fh^* \times \fg^*$ be a fixed
dynamical Lie algebroid corresponding to an $R : \fh^* \to 
\call (\fg^*,\fg)$ which satisfies the conditions
of Proposition \ref{pro:R} of the last section.
  In what follows, we shall formulate our results for the differentiable
category, but it will be clear that the results are also valid for 
the holomorphic category.

\begin{defi}
\label{def:3.1}
{\rm A Poisson manifold $(X,\pi_X)$ is said to admit a realization in
the dynamical Lie algebroid $T^*\fh^* \times \fg^*$ if there is a
Poisson map $\rho : X \to T\fh^* \times \fg$, where $T\fh^* \times
\fg$ is the dual vector bundle of $T^* \fh^* \times \fg^*$ equipped
with the Lie-Poisson structure.}
\end{defi}

\begin{defi}
{\rm Suppose a Poisson manifold $(X,\pi_X)$ admits a realization $\rho
: X \to T\fh^* \times \fg$ and $\HH \in C^{\infty} (X)$.  We say that
the Hamiltonian system $\dx = X_\HH (x)$ is realized in $T\fh^* \times
\fg$ by means of $\rho$ if there exists $K \in C^{\infty} ( T\fh^*
\times \fg)$ such that $\HH = \rho^*  K $.}
\end{defi}

In the following discussion, we shall work with a Poisson manifold
$(X, \pi_X   )$ together with a realization $\rho : X \to T\fh^* \times
\fg$.  Let $Pr_1 : T\fh^* \times \fg \to T \fh^*$, $Pr_2 : T\fh^* \times
\fg \to \fg$ be the  projection maps onto the first and second factor of
$T\fh^* \times \fg$ respectively and set
\begin{equation}
\label{eq:3.3}
	L = Pr_2 \circ \rho : X \to \fg;
\end{equation}
\begin{equation}
\label{eq:3.4}
	\tau = Pr_1 \circ \rho : X \to T \frakh^* .
\end{equation}
We also put
\begin{equation}
\label{eq:3.5}
	m = p \circ \tau : X \to  \frakh^* ,
\end{equation}
where $p : T \frakh^*  \to \frakh^*  $ is the bundle projection.  The next Proposition
shows how to compute the Poisson brackets of pullback of functions
in $Pr_2^* C^{\infty}(\fg)$ under the map $\rho$.  It is a direct
consequence of the canonical character of $\rho$ and the definition
of the Lie algebroid bracket on $T^*\fh^* \times \fg^*$.

\begin{pro}
\label{pro:3.6}
For all $f,g \in C^{\infty}(\fg )$, we have
\begin{eqnarray}
	& & \{ L^*f, \ L^*g \}_X (x)   \label{eq:LL} \\
	& = & \langle L(x), \ -ad_{R_{m(x)}(df(L(x)))}^{*} dg(L(x)) +
		ad_{R_{m(x)}(dg(L(x)))}^{*}  df(L(x))\rangle  \nonumber \\
	& & + \langle (\calx_{\tau(x)}R)(df(L(x))), \ dg(L(x))\rangle,\quad\quad
		\forall x \in X.  \nonumber
\end{eqnarray}
Here, and in the sequel, $df(L(x))$ and $dg(L(x))$  are 
considered as  elements in $\frakg^*$ for any fixed $x\in X$.
\end{pro}
\pf For any $\xi \in \fg^*$, we let $\ell_{\xi}$ denote the
corresponding linear function on $\fg$.  Then 
we have
\begin{eqnarray}
	& & \{ L^*f,\ L^*g \}_X (x)  \nonumber \\
	& = & \{ L^*\ell_{df(L(x))},  \ L^* \ell_{dg(L(x))}\}_X (x)\\
	& = & \{ Pr_2^*\ell_{df(L(x))}, \ Pr_2^*\ell_{dg(L(x))}\} ( \rho (x)) 
 \ \ \ (\mbox{since $\rho$ is a Poisson map}) \nonumber   \\
	& = & \langle [ df(L(x)),dg(L(x))] (m(x)), \  \rho(x)\rangle
\ \ \ (\mbox{by the definition of Lie-Poisson structure}) \nonumber   \\
	& = & \langle \tau(x),\  d \langle R(df(L(x))), \ dg(L(x))\rangle\ra
		  \nonumber \\
	& & + \langle L(x),\  - ad_{{R_{m(x)}(df(L(x)))}}^{*} dg(L(x)) +
		ad_{R_{m(x)}(dg(L(x)))}^* df(L(x))\rangle.   \qquad
		   \nonumber\\
& = & \langle L(x), \ -ad_{{R_{m(x)}(df(L(x)))}}^{*} dg(L(x)) +
                ad_{{R_{m(x)}(dg(L(x)))}}^*  df(L(x))\rangle  \nonumber \\
        & & + \langle \calx_{\tau(x)}R)(df(L(x))), \ dg(L(x))\rangle,\quad\quad
                \forall x \in X.  \nonumber
\end{eqnarray}
In the above computation, the quantities
  $df(L(x))$ and $ dg(L(x))$ are considered
as fixed elements in $\frakg^*$, the bracket
$[df(L(x)),dg(L(x))]$ is the Lie algebroid bracket when
both $df(L(x))$ and $ dg(L(x))$ are considered
as constant sections of $T^* \fh^* \times \fg^*$,
and in the  second from the last equality,
   $\langle R(df(L(x))),dg(L(x))\rangle$
is considered as a function on $\frakh^*$ with $x$ being fixed.\qed

\begin{rmk}
\label{rmk:3.1.}
 {\em If $R_{q}=r(q)^{\#}\in \call (\fg^* , \fg )$ 
as in Remark 2.1, then Equation 
(\ref{eq:LL}) is equivalent to the following
fundamental Poisson bracket relation:
$$\{L\ \stackrel{\ot}{,}\   L\} =[r^{12},\  L^{1}+ L^{2}]-\tau (x)r,$$
where $L^{1}=L\ot 1$ and $ L^{2}=1\ot L$. }
\end{rmk}

Let $I(\fg)$ be the collection of smooth ad-invariant functions on
$\fg$, i.e. $f \in I(\fg)$ iff 
 $ad_p^* df(p) = 0$ for all $p \in \fg$.
A natural collection of functions on $T\fh^* \times \fg$ is
$Pr_2^* I( \fg )$, the pullback of ad-invariant functions on $\fg$
by the projection map $Pr_2$.  As the reader will see, these
functions do not Poisson commute with respect to the Lie-Poisson
structure on $T\fh^* \times \fg$.  Thus our situation here is
quite different from that in standard classical $r$-matrix theory
for constant $r$-matrices.  We now examine the Hamiltonian systems
$\dx = X_{\HH}(x)$ on $X$ which can be realized in $T\fh^* \times \fg$
by means of $\rho$ with $\HH \in \rho^* (Pr_2^* I(\fg )) = L^{* } I( \fg )$.

\begin{pro}
\label{pro:3.8}
\begin{enumerate}
\item  If $\HH = L^* f$, where $f \in I(\fg)$, then under the flow
$\phi_t$ generated by the Hamiltonian $\HH$, we have the quasi-Lax type
equation:
\begin{eqnarray}
	\frac{dL(\phi_t)}{dt} & = & [R_{m(\phi_t)}(df(L(\phi_t))), \ 
		L(\phi_t)]  \nonumber \\
	& & - (\calx_{\tau(\phi_t)} R)(df(L(\phi_t))). \label{eq:3.9}  
\end{eqnarray}
\item  \ For all $f_1,f_2 \in I(\fg )$, we have
\begin{equation}
	\{ L^*f_1, L^*f_2\}_{X} (x) = \la  (\calx_{\tau(x)} R)(df_1 (L(x))), \ 
		df_2 (L(x)) \ra , \qquad \forall x \in X.
\end{equation}
\end{enumerate}
\end{pro} 
\pf 
(1) \  Let $\pi_X^{\#}:T^*X\lon TX$  be the  induced bundle map
of the Poisson tensor $\pi_X$ defined by $<\pi_X^{\#} \alpha  ,\ \beta>
=\pi_X (\alpha  ,\ \beta  ) , \ \forall \alpha,  \beta \in T^*X$.
From the invariance property of $f$ and Equation  (\ref{eq:LL}),
 we have
\begin{eqnarray}
	& & (T_x L\smalcirc \pi_X^{\#}(x) \smalcirc T_x^* L) [df(L(x))]   \nonumber \\
	& = & - [ R_{m(x)}(df(L(x))), \ L(x)] + (\calx_{\tau(x)} R)(df(L(x))),
	\nonumber
\end{eqnarray}
from which the assertion follows.

(2) \ This is obvious from Equation  (\ref{eq:LL}) and the invariance
property  of $f_1,f_2$.\qed

\begin{rmk} 
\label{rmk:3.2}
{\em  It is clear that the functions in $Pr_2^* I(\fg)$ do not Poisson
commute, for otherwise, it would contradict Proposition \ref{pro:3.8} (2).}
\end{rmk}

Proposition \ref{pro:3.8} (2) shows that there is an obstruction for 
$L^*I(\fg)$ to give a Poisson commuting family of functions.
A naive way to get rid of this obstruction is to restrict to 
the submanifold\\ $\tau^{-1}$ (zero section of $T\fh^*$).  It
is easy to see that $\tau^{-1}$ (zero section of $T\fh^*$)
 is a coisotropic submanifold of $X$ as the zero section
of $T\fh^*$ is a coisotropic submanifold of $T\fh^{*} $.  Thus one can
obtain a Poisson bracket on the quotient of $\tau^{-1}$ (zero
section of $T\fh^*$) by the characteristic foliation.  
Unfortunately, it is not necessary that  $\HH \in L^* I(\fg)$ 
or $L : X \to \fg$ will descend to the quotient space.  In the
following, we shall describe a situation where we can obtain 
integrable flows on a reduced phase space.
\newcounter{bean}
Let $\DD $ be a Lie subgroup of $G$ corresponding to
the Lie algebra $\frakh$.  We shall make the following 
assumptions:

\begin{list}
{A\arabic{bean}}{\usecounter{bean}
	\setlength{\rightmargin}{\leftmargin}}
\item $X$ is a Hamiltonian $\DD$-space with an equivariant momentum
	map $J: X \to \fh^*$,
\item the realization map $\rho : X\lon T\fh^* \times \frakg$
is equivariant, where $H$ acts on $T\fh^* \times \frakg$ by
adjoint action on the second factor.
\item there exists an $\DD$-equivariant map $g:X \to \DD$, where
$H$ acts on itself  by left translation, i.e.,
$$g(d \cdot x) = d \cdot g(x), \ \  d \in \DD, x \in X. $$
\end{list}

Suppose $\mu \in \fh^*$ is a regular value of $J$.
  Then, under the assumption
that $X_{\mu} = J^{-1} (\mu)/\DD$ is a smooth manifold, it
follows by Poisson reduction \cite{MR} that $X_{\mu}$ inherits a
unique Poisson structure $\{\cdot,\cdot\}_{X_{\mu}}$ satisfying
\begin{equation}
\label{eq:3.10}
	\pi^{*} \{\phi,\psi\}_{X_{\mu}}  = i^* \{\wphi,\wpsi\}_X  .
\end{equation}
Here, $i : J^{-1}(\mu) \to X$ is the inclusion map, $\pi : J^{-1}
(\mu) \to X_{\mu}$ is the canonical projection; $\phi, \psi
\in C^{\infty}(X_{\mu})$, and $\wphi,\wpsi$ are (locally defined)
smooth extensions of $\pi^* \phi,\pi^*\psi$ with differentials vanishing
on the tangent spaces of the $\DD$-orbits.  
It follows from Assumption A2 that 
$L:X \to \fg$ is $\DD$-equivariant, where the $\DD$-action 
        on $\fg$ is via the Ad-action.
Thus, if $\HH \in L^* I(\fg)$, 
it is clear that $\HH$ is $\DD$-invariant, so that $\HH$ descends
to  a function on $X_{\mu}$, i.e., there 
 exists a uniquely determined $\HH_{\mu}
\in C^{\infty}(X_{\mu})$ satisfying $\pi^* \HH_{\mu} = \HH \big|_{ J^{-1}
(\mu )}$.  However, as $L$ is only $\DD$-equivariant,
 therefore $L$ does not pass to the quotient 
and this is where Assumption A3 comes into play.
  Using the $\DD$-equivariant map $g$, we can define the gauge
transformation of $L$:
\begin{equation}
\label{eq:3.11}
	\wl : X \to \fg, \ \ \  x \mapsto Ad_{g(x)^{-1}} L(x)\,. 
\end{equation}
The following lemma is obvious.

\begin{lem}
\label{lem:3.12}
\begin{enumerate}
\item $\wl$ is $\DD$-invariant.
\item If $\HH \in L^* I(\fg)$, say, $\HH = L^* f$, then also $\HH = \wl^* f$.
\end{enumerate}
\end{lem}

It follows from this lemma that there exists a uniquely determined
map $L_{\mu} : X_{\mu} \to \fg$ such that
\begin{equation}
\label{eq:3.13}
	L_{\mu} \smalcirc \pi =  \wl \big|_{ J^{-1}(\mu )}\,.
\end{equation}
In particular, if $\HH = L^* f$, where $f \in I(\fg)$, then $\HH$
descends to a function $\HH_{\mu}$ on $X_{\mu}$ such that
\begin{equation}
\label{eq:3.14}
	\HH_{\mu} = L_{\mu}^* f.
\end{equation}
In other words, the functions in  $L^* I(\fg)\;\big|_{J^{-1}(\mu)}$
descends into functions in  $L_{\mu}^*I({\fg}) \subset C^{\infty}
(X_{\mu})$.

The following lemma is straightforward from the definition of
$\wl$ in Equation (\ref{eq:3.11}).

\begin{lem}
\label{lem:3.15}
$T_x \wl = Ad_{g(x)^{-1}} \smalcirc  T_x L + ad_{\wl(x)}  \smalcirc
 l_{g(x)^{-1}* } \smalcirc
	T_x g$, $\forall x \in X$, where both sides are
considered as linear maps from $T_{x}X $ to $\frakg$,
and $l_{g(x)^{-1} }$ is left translation by $g(x)^{-1}\in H$. 
\end{lem}


We now make an additional assumption. \\\\

A4\  $\calx_v  R = 0$, $\forall v\in \tau (J^{-1}(\mu ))$.\\\\

\begin{pro}
\label{lem:3.16}
\begin{equation}
T_x \wl  \smalcirc \pi_X^{\#}(x)  \smalcirc T_x^*  \wl = ad_{\wl(x)} 
\smalcirc  \wr(x) + \wr(x) \smalcirc 
	ad_{\wl(x)}^*\,,  \ \forall x \in J^{-1}(\mu),
\end{equation}
 where both sides of the equation are considered
as  linear maps from $\frakg^*$ to $\frakg$,  where
$\wr : J^{-1}(\mu) \to \call (\fg^*, \ \fg)$ is given by
\begin{eqnarray}
	\wr (x) & = & Ad_{g(x)^{-1}} \smalcirc (R \smalcirc m)(x)
\smalcirc  Ad_{g(x)^{-1}}^*
	+ T_x L  \smalcirc \pi_X^{\#}(x) \smalcirc  T_x^* g \smalcirc
l_{g(x)^{-1}}^{*}
\label{eq:3.17} \\
	& + & \frac12\; ad_{\wl(x)}  \smalcirc  l_{g(x)^{-1}*} \smalcirc 
T_x g  \smalcirc \pi_X^{\#}(x) \smalcirc  T_x^* g \smalcirc
l_{g(x)^{-1}}^{*}, \ \ \forall x\in J^{-1}(\mu ). \nonumber
\end{eqnarray}
Moreover, $\wr $ is $\DD$-invariant. Here, as well as in the sequel,
$Ad^*$  denotes the dual map of $Ad$ defined by:
$\la Ad^*_d \xi, \ X\ra =\la \xi, \   Ad_d X\ra$, $\forall \xi\in
\fg^*$ and $X\in \fg$.
\end{pro}
\pf   Apply Lemma \ref{lem:3.15} and Proposition \ref{pro:3.6},
 together with A4, the
expression for $T_x \wl \smalcirc \pi_X^{\#}(x) \smalcirc  T_x^* \wl$ follows.
On the other hand, it follows from
Assumption  A2  ($\rho$ is $H$-equivariant) that $m (d\cdot x)=
m(x)$, $\forall x \in X$, $d \in \DD$. Hence,
\begin{equation}
\label{eq:Rm}
(R \smalcirc m) (d \cdot x)=(R \smalcirc m)(x)=
Ad_d \smalcirc  (R\smalcirc m)(x) \smalcirc Ad_{d}^* , 
\end{equation}
since $(R \smalcirc m)(x)\in \call (\fg^* , \fg)^H$ according
to Proposition  \ref{pro:R}.
 Thus, the assertion that $\wr$ is $\DD$-invariant is a consequence
of  Equation (\ref{eq:Rm}),
the equivariance property of the maps $L$ and $g$.
  We shall omit the straightforward calculations.   \qed

From the  $H$-invariance of $\wr$, it follows that there exists
$R_{\mu}: X_{\mu} \to \call (\fg^*, \ \fg)$ such that
\begin{equation}
\label{eq:3.18}
	R_{\mu} \smalcirc \pi = \wr\,.
\end{equation}

We now come to the main result of the section.

\begin{thm}
\label{thm:3.19}
Let $(X, \ \pi_X)$ be a Poisson manifold with a realization $\rho : X \to
T\fh^* \times \fg$ which satisfies A1--A4. 
  Then, under the assumption that $X_{\mu} = J^{-1}(\mu)/\DD$ 
is a smooth manifold, there
exists a unique Poisson structure $\{\cdot, \ \cdot\}_{X_{\mu}}$ on
$X_{\mu}$ satisfying Equation (\ref{eq:3.10})
 and a map $L_{\mu} : X_{\mu} \to \fg$
satisfying  Equation (\ref{eq:3.13}) such that
\begin{enumerate}
\item  $\forall  f_1, \ f_2 \in C^{\infty} (\fg)$,
\begin{eqnarray}
	& & \{L_{\mu}^* f_1, L_{\mu}^* f_2\}_{X_{\mu}}(\xx) \nonumber \\
	& = & - \langle L_{\mu}(\xx), \ ad^{*}_{R_{\mu}^{*}(\xx)
		(df_2(L_{\mu}(\xx)))} df_1(L_{\mu}(\xx))  \nonumber \\
	& & \qquad + ad_{R_{\mu}(\xx)(df_1 (L_{\mu}(\xx)))}^*
		df_2(L_{\mu}(\xx))\rangle\,, \ \ \forall \xx\in X_{\mu}; 
\label{eq:LLu}
\end{eqnarray}
\item 
 Functions in $L_{\mu}^* I(g)$ Poisson commute in $(X_{\mu},
	\{\cdot,\cdot\}_{X_{\mu}})$;
\item  If $\HH_{\mu} = L_{\mu}^* f$, $f \in I(\fg)$, then under the
	flow generated by $\HH_{\mu}$, we have
\begin{equation}
	\frac{d\,L_{\mu}}{dt} = - [(R_{\mu})^{*}(df(L_{\mu})), \ 
	L_{\mu}].
\end{equation}
\end{enumerate}
\end{thm}
\pf
\noindent(1). \ Let $\xx=\pi ({x})\in X_{\mu }$
for some $x\in J^{-1}(\mu )$.
 From Equations (\ref{eq:3.10}), (\ref{eq:3.13}), 
(\ref{eq:3.18}) and  Proposition \ref{lem:3.16}, we have
\begin{eqnarray}
	& & \{L_{\mu}^* f_1, \ L_{\mu}^* f_2\}_{X_{\mu}} (\xx)  \nonumber \\
	& = & \{\wl^* f_1,  \ \wl^* f_2\}_X ({x})  \nonumber  \\
	& = & \langle (ad_{\wl({x})} \smalcirc \wr({x})
 + \wr({x})   \smalcirc ad_{\wl({x})}^*)
		(df_1(\wl({x}))), \ \  df_2(\wl({x}))\rangle  \nonumber
	\\
	& = & - \langle ad_{R_{\mu}(\xx)(df_1(L_{\mu}(\xx)))} L_{\mu}
	(\xx), \ \  df_2 (L_{\mu}(\xx))\rangle  \nonumber
	\\
	&  & - \langle df_1(L_{\mu}(\xx)) , \ \ 
ad_{R_{\mu}^{*}(\xx)(df_2(L_{\mu}(\xx)))} L_{\mu} (\xx)\rangle, \nonumber
\end{eqnarray}
from which the  assertion follows.

\noindent(2). \ This is obvious from (1).

\noindent(3). \ If $\pi_{X_{\mu}}^{\#}$ denotes  the induced
 bundle map $T^* X_{\mu} \lon TX_{\mu}$
of the Poisson tensor  on $X_{\mu}$, it follows from (1) and the invariance
property of $f$ that $(TL_{\mu} \smalcirc \pi_{X_{\mu}}^{\#} \smalcirc T^* L_{\mu})
(df(L_{\mu})) = [(R_{\mu})^* (df(L_{\mu})), \ L_{\mu}]$.  Hence the
assertion is immediate.   \qed

\begin{rmk}
{\em  If $R =r^{\#}: \frakh^* \lon \call (\frakg^*, \ \frakg)$ 
for a classical  dynamical r-matrix as in Remark \ref{rmk:2.1},
 then Equation (\ref{eq:LLu}) is equivalent to the following
relation:
\begin{equation}
 \{L_{\mu} \ \stackrel{\ot}{,} \  \ L_{\mu}\} (\xx)
=[\tilde{r}^{12}(\xx), L_{\mu}^{1} (\xx)]
-[\tilde{r}^{21}(\xx), L_{\mu}^{2}(\xx)]\ \ \ \forall \xx \in X_{\mu},
\end{equation}
where 
\begin{equation}
\label{eq:tr}
\tilde{r}(\xx)=Ad_{g({x})^{-1}}(r(m( {x}))-\{g^{1}, L^{2}\}
({x}){g^{1}}^{-1}+\half [u^{12 }({x}), L^{2}({x})]).
\end{equation}
Here,  $ {x}\in J^{-1}(\mu )$  is such that
$\xx=\pi ({x})$,
 $u^{12}=(g_{*}\pi_{X})g^{-1}\in C^{\infty}(X, \wedge^{2}\frakg )$,
and $\{g^{1}, L^{2}\}{g^{1}}^{-1}\stackrel{def}{=} \half \sum
((g_{*}X_{i})g^{-1}\ot L_{*}Y_{i} -(g_{*}Y_{i})g^{-1}\ot L_{*}X_{i})$
as a map from $X$ to $\frakg\ot \frakg$, where
$X_{i}, Y_{i}\in {\frak{X} } (X)$ are $\DD$-invariant vector
fields such that $\pi_{X} = \sum  X_{i}\wedge Y_{i}= \half
\sum (X_{i}\ot Y_{i}-Y_{i}\ot X_{i})  $, $X_{i}, Y_{i}\in {\frak X} (X)$.

We remark that fundamental Poisson bracket relations
of this nature, in which the r-matrix can
depend on phase space variables, was first considered
in \cite{BV}. }
\end{rmk}

\section{Spin   Calogero-Moser systems}

Let $\frakg$ be a Lie algebra over $\complex$ with a
non-degenerate ad-invariant bilinear form $(\cdot , \cdot )$
 and $\frakh\subset \frakg$
 a non-degenerate  (i.e., the restriction of $(\cdot , \cdot )$
to $\frakh$ is non-degenerate) Abelian Lie     subalgebra.
 By definition (see Remark \ref{rmk:4.1} below),   a  classical
dynamical $r$-matrix with  spectral parameter associated with the
pair $\frakh\subset \frakg$
is  a meromorphic  map
$r \, : \, \h^* \,\times \,\complex \, \to \g \otimes \g $
having a  simple pole at $z=0$ and satisfying the following conditions:

\begin{enumerate}
\item the zero weight condition:
\begin{equation}
\label{eq:4.1}
[h\otimes 1 + 1\otimes h\, , \, r(q,  z)]=0,
\end{equation} 
for all $h \in \h$ and  all $(q , z) \in \h^* \times  \complex $ 
except for the poles of $r$;
\item the generalized unitarity condition:
\begin{equation}
\label{eq:4.2}
r^{12}(q,z)\, + \,r^{21}(q ,-z) \,= \, 0,
\end{equation}
for all $(q , z) \in \h^* \times  \complex $ 
except for the poles of $r$;
\item the residue condition:
\begin{equation}
\label{eq:res}
 {\mbox{Res}}_{z=0} \, r(q, z)\, = \, \Omega,
\end{equation}
where  $\Omega \in (S^{2}\frakg )^{\frakg}$ is the  Casimir 
element corresponding
to the bilinear form $(\cdot, \cdot )$;
\item the classical dynamical Yang-Baxter equation (CDYBE):
\begin{equation}
\label{eq:cdybe}
\Alt  (d_\h r) \, + \, [r^{12}(q, z_{1,2}), r^{13}(q,z_{1,3})]\,+ \,
[r^{12}(q,z_{1,2}), r^{23}(q, z_{2,3})]\,
+\,[r^{13}(q, z_{1,3}), r^{23}(q,z_{2,3})] \, = \, 0 \, ,
\end{equation}
where $z_{i,j}=z_i-z_j$.

In Equation (\ref{eq:cdybe}), the differential
of the $r$-matrix is considered with respect to the
$\h^*$-variables:
$$d_\h r :\h^* \times  \complex  \lon \frakg\ot \frakg\ot \frakg,
\ \ \ \ (q, z)\lon  \sum _i \, h_i^{(1)} \ot {\partial r^{23}
\over \partial q_i}(q,z), $$
and the term  $\Alt  (d_\h r)$ is a shorthand for
the following symmetrization  of $d_\h r$:
\begin{equation}
\label{eq:4.5}
\Alt (d_\h r) \, = \, \sum _i \, h_i^{(1)} \ot {\partial r^{23}
\over \partial q_i}(q,z_{2,3})\, + \, \sum _i \, h_i^{(2)}  \ot{\partial r^{31}
\over \partial q_i}(q,z_{3,1})\, + \, \sum _i \, h_i^{(3)}  \ot{\partial r^{12 }
\over \partial q_i}(q, z_{1,2})\, ,
\end{equation}
where $(h_{1}, \cdots , h_{N})$ is a basis of $\frakh$,
and $(q_{1}, \cdots , q_{N})$  its corresponding coordinate
system on $\frakh^*$. 
\end{enumerate}

We call the variable  $z$  in $r (q, z)$ the spectral parameter.

By $L\frakg$, we denote the Lie algebra of Laurent series
$X=\sum_{n=-T}^{\infty}X_{n}z^n$ with coefficients
in $\frakg$, which are convergent in some annulus
$A_{c}=\{z\in \complex |0<|z|<c \}$ (which may depend on the series).
The Lie bracket in $L\frakg$ is the pointwise  bracket.
In a similar fashion, we can define the restricted
 dual $L\frakg^*$.
Using the bilinear  form on $\frakg$, we can define
a non-degenerate invariant bilinear form on $L\frakg$ by
\begin{equation}
\label{eq:4.18}
(X , Y )= {\mbox{Res}}_{z=0} (X (z),\  Y (z )), \ \ \ \forall X, Y\in L\frakg.
\end{equation}
On the other hand, the pairing between $L\frakg^*$
and $L\frakg$ is given by
\begin{equation}
\label{eq:4.19}
\la \xi , X \ra = {\mbox{Res}}_{z=0} \la  \xi  (z),\  X (z ) \ra 
, \ \ \ \forall \xi \in
L\frakg^* , \ X\in  L\frakg.
\end{equation}
 

Associated with each dynamical $r$-matrix $r$ with spectral
parameter
 is  an operator $R: \frakh^* \lon \call (L\frakg^*  , L\frakg)^{H}$,
which we use to define a Lie algebroid structure
on $T^* \frakh^* \times L\frakg^* $ according to  the recipe
in Section 2. We now proceed with the construction
of $R$. Let
\begin{equation}
\label{eq:4.7}
r(q, z)=\frac{\Omega}{z}+\sum_{k=0}^{\infty}
t_{k}(q)z^{k}, \ \ \ \ \ \ z\in A_{c(r)}
\end{equation}
be the Laurent expansion of $r(q, \cdot )$ about $z=0$,
where $c(r)$ denotes the radius of convergence of the series.
Assume furthermore that \\\\

5.\ \ $c(r)$ is independent of $q$ (which we will always
assume in the sequel when talking about a dynamical r-matrix with spectral
parameter).\\\\\\

\begin{rmk}
\label{rmk:4.1}
{\em The original definition of classical dynamical r-matrices
with spectral parameter is for simple Lie algebras \cite{EV}.
In the above, we have modified this definition
by putting in the extra assumptions. Namely, the pole of 
$r(q, \cdot )$ at $z=0$ is simple and the number $c(r)$ is
independent of $q$. For simple Lie algebras, these
additional assumptions are not necessary as they follow from the
solution of the classification problem \cite{EV}.}
\end{rmk}

For any $\xi \in L\frakg^* $, denote by $A_{c(\xi )}$ the
largest annulus on which the Laurent series converges
and let $c_{0}(r, \xi )=\half min (c(r), c(\xi ) )$.
If $q\in \frakh^*$ is not a pole of $r(\cdot , z)$, we set
\begin{equation}
\label{eq:rR}
(R_{q} \xi)(z)=p.v. \frac{1}{2\pi i} \oint_{C} \la r(q, w-z),  \ \xi (w)
\otimes 1 \ra  dw,
\ \ \ \ \ \forall z\in A_{c_{0} (r, \xi )}, \ \xi\in L\frakg^*, 
\end{equation}
where $C$ is the circle centered at 0 of radius $|z|$
 with positive orientation,
and $p.v.$ denotes the principal value of the  improper integral.

\begin{lem}
\label{lem:4.8}
$R_{q} \xi$ is well-defined on $A_{c_{0}(r, \xi )}$, i.e.,
 the principal value of the  improper integral in 
Equation (\ref{eq:rR}) exists.
\end{lem}
\pf Consider a circle $K$ centered at $z\in A_{c_{0}(r, \xi )}$ with a
small radius $\epsilon$ such that $K$
intersects $C$ at exactly two points $z'$ and $z''$.
We denote by $C_{\epsilon}$ the circular arc $z' z''$ and
by $K'$ the portion of $K$ which lies
to the left of $C_{\epsilon}$ with orientation
given by the clockwise direction.
By definition,
\be
&&p.v. \frac{1}{2\pi i} \oint_{C} \la r(q, w-z), \  \xi (w)\otimes 1 \ra dw\\
&=& \frac{1}{2\pi i}
\lim_{\et}\int_{C-C_{\epsilon}} \la r(q, w-z), \ \xi (w)\otimes 1 \ra dw.
\ee

We have
\be
&&\int_{C-C_{\epsilon}} \la r(q, w-z), \ \xi (w)\otimes 1 \ra dw\\
&=&\int_{C-C_{\epsilon}} \la r(q, w-z),  \ (\xi (w)- \xi (z)) \ot 1 \ra dw
+ \int_{C-C_{\epsilon}} \la r(q, w-z),  \ \xi (z) \ot 1 \ra dw.
\ee

Since $\xi $ is analytic at $z$, it follows from the
residue condition:  (\ref{eq:res}) that

\be
&&\lim_{\et}
\int_{C-C_{\epsilon}} \la r(q, w-z), \  (\xi (w)- \xi (z)) \ot 1 \ra dw\\
&=&\int_{C}\la r(q, w-z), \  (\xi (w)- \xi (z)) \ot 1 \ra dw.
\ee

On the other hand,
\be
&&\int_{C-C_{\epsilon}} \la r(q, w-z), \ \xi (z)\otimes 1 \ra dw\\
&=&(\int_{C-C_{\epsilon}+K'}-\int_{K'}) \la r(q, w-z),  \ \xi (z)\otimes 1 \ra dw\\
&=&-\int_{K'}\la  r(q, w-z), \  \xi (z)\otimes 1 \ra dw,
\ee
because $\la r(q, w-z),  \ \xi (z)\otimes 1 \ra $, as a 
function of $w$, is analytic
in the interior of $(C-C_{\epsilon})+K'$. Now

\be
&&-\int_{K'}\la r(q, w-z),\  \xi (z)\otimes 1 \ra dw\\
&=&-\int_{K'} \la \frac{\Omega}{w-z}, \  \xi (z)\ot 1 \ra dw
-\int_{K'} \la \sum_{k=1}^{\infty}t_{k}(q)(w-z)^{k},  \   \xi (z)\ot 1 \ra dw\\
&=&-(I\xi ) (z)(\log \left|\frac{z''-z}{z'-z}\right|
+iVar_{K'}Arg(w-z))\\
&&-\int_{K'}\la \sum_{k=1}^{\infty} t_{k}(q)(w-z)^{k},   \  \xi (z)\ot 1 \ra dw\\
&\stackrel{\et}{\lon}& \pi i (I\xi ) (z),
\ee
where $I : L\frakg^* \lon L\frakg$ is the linear  isomorphism 
induced by the bilinear form $(\cdot , \cdot )$ as defined
by Equation (\ref{eq:4.18}). 

Consequently, the principal value of the improper integral in
Equation (\ref{eq:rR}) exists. \qed 

Indeed, from the proof of the above lemma, we obtain
the formula
\begin{equation}
\label{eq:4.9}
(R_{q}\xi )(z)=\half (I\xi ) (z) +
\frac{1}{2\pi i} \oint_{C}\la r(q, w-z), \ (\xi (w)- \xi (z)) \ot 1 \ra dw,
\ \ \ \ \forall \xi (z) \in L\frakg^* , \  z\in  A_{c_{0}(r, \xi )},
\end{equation}
which shows that $R_{q}\xi $ is analytic in the annulus $A_{c_{0}(r, \xi )}$.
We can therefore extend $R_{q}\xi$ to other possible values
of $z$ by analytic continuation. In this case, we can do it
explicitly using the following

\begin{pro}
\label{pro:4.10}
For $z\in  A_{c_{0} (r, \xi )}$, we have the formula
\begin{equation}
\label{eq:4.11}
(R_{q}\xi )(z)=\half (I\xi ) (z) +
\sum_{k\geq 0}\frac{1}{k!}\la \frac{\partial^{k}r}{\partial z^{k}}
(q, -z ), \ \xi_{-(k+1)}\ot 1 \ra .
\end{equation}
Hence we can analytically continue $R_{q}\xi $ to $A_{c(r, \xi )}$ by
using this formula, where $c(r, \xi )=min (c(r), c(\xi ))$.
\end{pro}
\pf Let $C$ be the circle centered at $0$ of radius $|z |$ with positive
orientation, and introduce the   map
$$\Phi (\lambda ) =\frac{1}{2\pi i} \oint_{C}\la r(q, w-\lambda ),\ 
 \xi (w) \ot 1 \ra  dw ,
\ \ \ \ \ \lambda \in A_{c_{0} (r, \xi )}-C .$$
If $\lambda $ is on the $+$-side of $C$ (i.e. the interior of
$C$), we have 
$$ \Phi (\lambda ) = \frac{1}{2\pi i} \oint_{C} \la r(q, w- \lambda ),  \ 
(\xi (w)- \xi ( \lambda )) \ot 1 \ra dw
+(I\xi ) ( \lambda ). $$

Therefore, the boundary value
\be
&&\Phi^{+} (z) \\
&=&\lim_{ \stackrel{\lambda \to z}{ \lambda \in + \mbox{side}} } \Phi (\lambda )\\
&=&\frac{1}{2\pi i} \oint_{C}\la r(q, w-z), \  (\xi (w)- \xi (z)) \ot 1 \ra dw
+(I\xi ) ( z), \ \ \ z\in C.
\ee

On comparing this equation with Equation (\ref{eq:4.9}), we obtain
\begin{equation}
\label{eq:4.12}
\Phi^{+}(z) =\half (I\xi ) ( z) + (R_{q}\xi )(z).
\end{equation}
But for $\lambda $ on the $+$-side of $C$, the integrand 
of $\Phi (\lambda )$  has poles at $w=\lambda $ and $w=0$
in the interior of $C$. Consequently, by
the residue theorem,
\be
&&\Phi ( \lambda )\\
&=& {\mbox{Res}}_{w= \lambda } \la r(q, w-\lambda  ), \xi (w) \ot 1 \ra 
+{\mbox{Res}}_{w=0} \la r(q, w-\lambda  ), \xi (w) \ot 1 \ra \\
&=& (I\xi ) (\lambda )+\sum_{k\geq 0} 
\frac{1}{k!} \la  \frac{\partial^{k}r}{\partial \lambda^{k}} (q, -\lambda ), \ 
 \xi_{-(k+1)}\ot 1  \ra .
\ee
 Equation (\ref{eq:4.11}) now follows from  letting
$\lambda \to z$ and Equation (\ref{eq:4.12}). \qed

We now define the operator $R: \ \fh^* \lon \call (L\fg^* , L\fg )$ by
\begin{equation}
\label{eq:Rq}
R(q)\xi =R_{q}\xi
\end{equation}
for all $q\in \fh^*$ which is not a pole of $r$ and
for all $\xi \in L\fg^*$. We shall use
Equation (\ref{eq:4.11}) to compute $R$ from $r$.

If $q\in \frakh^*$ is not a pole of $r(\cdot , z)$, we
define $r_{-}^{\# }(q): \frakg^{*} \lon L \frakg$ by
\begin{equation}
\la (r_{-}^{\# }(q)\xi )(z), \  \eta \ra =
\la r(q, z), \  \eta \ot  \xi \ra , \ \ \ \forall \xi , \eta \in \frakg^* .
\end{equation}
 From the generalized unitarity condition,
it is easy to check that $R_q^* =-R_q$. We now 
examine the consequences of the zero weight condition
and the classical dynamical Yang-Baxter 
equation  which are basic in our theory.

\begin{lem}
\label{lem:3.13a}
 Let $r: \frakh^* \times \complex \lon \frakg\ot \frakg$ be
a classical dynamical $r$-matrix with spectral parameter.
Then we have
\begin{enumerate}
\item $\la r(q, z), \  Ad_{x^{-1}}^* \xi \ot 1 \ra =Ad_{x} \la r(q, z), \  \xi \ot 1 \ra , \ \ \ x\in H, \ \ \xi\in \frakg^* $;
\item $  (r_{-}^{\# }(q))^{*} (ad_{h}^{*} \xi )=ad_{h} [(r_{-}^{\# }(q))^{*} 
\xi ], \ \ \ \ h\in \frakh, \ \ \xi\in L\frakg^* $;
\item $j^{*} (I^{-1} ([R_{q}\xi ,  I\eta ]+[ I\xi , R_{q} \eta ]))=0, \ \ \ 
\forall \xi, \eta \in L\frakg^{*} $, 
where $j:\frakh \lon L\frakg$ is
the natural inclusion, and $j^* : L\frakg^* \lon \frakh^* $ is
the dual map.
\end{enumerate}
\end{lem}
\pf (1) The relation is a simple consequence of the global
version of the zero weight condition.

(2) From the zero weight condition, it follows that for all
$h\in \frakh$, $\xi \in L\frakg^* $, $\eta \in \frakg^* $, we have
\begin{eqnarray}
0&=&\la r(q, z), \ ad_{h}^{*} \xi (z)\ot \eta + \xi (z)\ot ad_{h}^{*} \eta \ra
\nonumber \\
&=&\la (r_{-}^{\# }(q)\eta )(z), \  ad_{h}^{*} \xi (z) \ra 
+\la (r_{-}^{\# }(q)ad_{h}^{*} \eta )(z), \  \xi (z) \ra . \label{eq:rxe}
\end{eqnarray}
If $C\subset A_{c(r, \xi )}$ is a circle centered at $0$
with positive orientation and we integrate the above relation
with respect to $z$ over $C$, the result is
$$\la r_{-}^{\# }(q)\eta ,\  ad_{h}^{*} \xi \ra +
\la r_{-}^{\# }(q)ad_{h}^{*} \eta , \ \xi \ra  =0, $$
by the definition of the pairing in Equation (\ref{eq:4.19}). As the 
above equality holds for all $\eta \in \frakg$, the assertion
follows.

(3) In Equation (\ref{eq:rxe}),  replace $r(q, z)$ by $r(q, z-w)$
and $\eta$ by $\eta (w)$, we have
\be
0&=&\la r(q, z-w), \ ad_{h}^* \xi (z)\ot \eta (w) +\xi (z)\ot ad_{h}^* 
\eta (w) \ra ,
\ee
$\forall h\in \frakh$, $\xi  , \  \eta \in L\frakg^* $. Let
$C\subset A_{c_{0}(r, \xi )} \cap A_{c (\eta )}$ be
a circle centered at $0$ with positive orientation.
For $w\in C$, take the principal value of the 
integral of the above expression with respect to
$z$ over $C$, we have
$$0=\la (R_{q}ad_{h}^{*} \xi )(w), \  \eta (w) \ra +
\la (R_{q}\xi )(w),  \  ad_{h}^{*}\eta (w) \ra .$$
Then an  integration  with respect to $w$ over $C$ yields
\be
0&=& \la R_{q}ad_{h}^{*} \xi , \eta \ra + \la R_{q}\xi  , ad_{h}^{*} \eta \ra\\
&=&-(j (h),  \ [R_{q}\xi ,  I\eta ]+[I\xi , R_{q} \eta ])\\
&=&-\la j (h),  \ I^{-1}([R_{q}\xi ,  I\eta ]+[I\xi , R_{q} \eta ])\ra
\ee
Therefore, $ j^{*}( I^{-1}([R_{q}\xi , I \eta ]+[I\xi , R_{q} \eta ]))=0$. \qed

To prepare for the proof of the next proposition, we first note
by a direct calculation that

\begin{eqnarray}
&&\la [r^{12}(q, z-w), r^{13}(q, z-v)], \ \xi (z)\ot \eta (w) \ot \zeta (v) \ra
\nonumber \\
&=&\la \xi (z) , \ [\la r(q, z-w), 1\ot \eta (w) \ra, \ \la r(q, z-v), 
 1\ot \zeta (v) \ra ] \ra \label{eq:3.14a} ,
\end{eqnarray}
\begin{eqnarray}
&&\la [r^{12}(q, z-w), r^{23}(q, w-v)], \ \xi (z)\ot \eta (w) \ot \zeta (v) \ra
\nonumber  \\
&=&\la \eta (w), \ [\la r(q, z-w), \xi (z)\ot 1\ra , \ \la r(q, w-v),
  1\ot \zeta (v) \ra ] \ra, \ \ \ \ \mbox{and} \label{eq:3.14b}
\end{eqnarray}
\begin{eqnarray}
&&\la [r^{13}(q, z-v), r^{23}(q, w-v )], \ \xi (z)\ot \eta (w) \ot \zeta (v) \ra
\nonumber  \\
&=&\la \zeta (v), \ [ \la r(q, z-w), \xi (z)\ot 1 \ra , \  \la r(q, w-v),
  \eta (w)\ot 1 \ra ] \ra , \label{eq:3.14c} 
\end{eqnarray}
where $\xi, \eta, \zeta \in L\frakg^* $.

Let $\pi_{\frakh}$ be the  projection operator onto  $\frakh$ relative
to the decomposition $\frakg =\frakh \oplus \frakm$, where
$\frakm$ is the orthogonal complement of $\frakh$.
  For any $\xi , \eta  \in L\frakg^*$,  let
$\pi_{\frakh} (I\xi) (z )=\sum_i \xi_{i}(z)h_{i}$ and
$\pi_{\frakh} (I\eta) (z )=\sum_i \eta_{i}(z)h_{i}$. We have the 
following relations corresponding  to the terms
in $\Alt (d_{\frakh} r)$:

\begin{eqnarray}
&&\la \sum_{i} h_{i}^{(1)}\ot \frac{\partial r^{23}}{\partial q_{i}}(q,  w-v),
\ \xi (z) \ot \eta (w )\ot \zeta (v)  \ra \nonumber\\
&=&\la \zeta (v), \ \sum_{i} \xi_{i} (z)\frac{\partial}{\partial q_{i}}
\la r(q,  w-v), \eta (w ) \ot 1  \ra \ra ; \label{eq:3.15a}\\
&&\la \sum_{i} h_{i}^{(2)}\ot \frac{\partial r^{31}}{\partial q_{i}}(q,  v-z),
\ \xi (z) \ot \eta (w )\ot \zeta (v) \ra  \nonumber \\
&=&-\la \zeta (v), \ \sum_{i} \eta_{i} (w)\frac{\partial}{\partial q_{i}}
\la r(q,  z-v), \xi (z) \ot 1 \ra \ra , \ \ \ \mbox{and} \label{eq:3.15b}\\
&&\la \sum_{i} h_{i}^{(3)}\ot \frac{\partial r^{12}}{\partial q_{i}}(q,  z-w),
\ \xi (z) \ot \eta (w )\ot \zeta (v) \ra   \nonumber\\
&=&\la \zeta (v), \ \sum_{i} h_{i} \frac{\partial}{\partial q_{i}}
\la r(q,  z-w), \xi (z) \ot \eta (w ) \ra \ra .  \label{eq:3.15c}
\end{eqnarray}

\begin{pro}
\label{pro:Rmdybe}
For each $q\in \frakh^*$ which is not a pole  of $r(\cdot , z)$,
  the operator $R_{q}$ 
is  in $\call (L\frakg^* , L\frakg)^{H}$ and 
satisfies the mDYBE (Equation (\ref{eq:mdybe})) with $c=-\frac{1}{4}$.
\end{pro}
\pf Let $0<c<1$ and let $ C\subset A_{c_{0}(r, \xi )}\cap A_{c(\eta )}
\cap  A_{c(\zeta  )}$
be a  circle centered at zero with 
positive  orientation. From Equation (\ref{eq:3.14a}), we have
\be
&&(\frac{1}{2\pi i})^{3} \oint_{C}  \lim_{\ct 1-} \oint_{C} p.v.
\oint_{C} \la  [r^{12}(q, z-w), r^{13}(q, z-cv )], \  \xi (z) \ot \eta (w )
\ot \zeta (cv)  \ra  dz dw dv\\
&=&(\frac{1}{2\pi i})^{3} \oint_{C}  \lim_{\ct 1-} \oint_{C} p.v.
\oint_{C}  \la \xi (z), \  [ \la r(q, z-w), 1\ot \eta (w) \ra ,
 \  \la  r(q, z-cv ), 1\ot \zeta (cv) \ra ] \ra  dw dz  dv\\
&=& (\frac{1}{2\pi i})^{2} \oint_{C}  \lim_{\ct 1-} \oint_{C} 
\la \xi (z), \ [(R_{q}\eta )(z), \  \la   r(q, cv-z ), \ \zeta (cv)\ot 1 \ra 
] \ra dz dv\\
&& \ \ \ (\mbox{by the generalized unitarity condition})\\
&=& -(\frac{1}{2\pi i})^{2} \oint_{C}\lim_{\ct 1-} \oint_{C}
\la \zeta (cv), \ \la r(q, z-cv ),  \ I^{-1}[I \xi , R_{q}\eta ](z )\ot 1 \ra \ra 
dz dv\\ 
&& \ \ \ (\mbox{by the ad-invariance of $(\cdot , \cdot )$
and the generalized unitarity condition})\\
&=&-\la \zeta , (R_{q}+ \half  I)(I^{-1} [I\xi , R_{q}\eta ] ) \ra 
 \ \ \ (\mbox{by Equation (\ref{eq:4.12})}). 
\ee
Note that we have interchanged the order of integration in going
from the first line to the second line of the above calculation.
This fact can be easily verified  and we leave the details 
to the reader. In what follows, it is not necessary
to interchange the order of integrations. Indeed, a similar
manipulation using Equation (\ref{eq:3.14b}) shows that
\be
&&(\frac{1}{2\pi i})^{3} \oint_{C}  \lim_{\ct 1-} \oint_{C} p.v.
\oint_{C}  \la [r^{12}(q, z-w), r^{23}(q, w-cv )], \  \xi (z) \ot \eta (w )
\ot \zeta (cv) \ra  dz dw dv\\
&=& \la \zeta , (R_{q}+\half I) (I^{-1} [ I\eta , R_{q}\xi ]) \ra .
\ee
Meanwhile, by using Equations (\ref{eq:3.14c}) and (\ref{eq:4.12}),
we find
\be
&&(\frac{1}{2\pi i})^{3} \oint_{C}  \lim_{\ct 1-} \oint_{C} 
\oint_{C}  \la [r^{13}(q, z-cv), r^{23}(q, w-cv )], \  \xi (z) \ot \eta (w )
\ot \zeta (cv) \ra  dz dw dv\\
&=& \la \zeta ,  [(R_{q}+\half I )  \xi , \ (R_{q} +\half I )\eta ] \ra .
\ee
On the other hand,  from Equation (\ref{eq:3.15a}), we have
\be
&&(\frac{1}{2\pi i})^{3} \oint_{C}  \lim_{\ct 1-} \oint_{C} 
\oint_{C}  \la \sum _i \, h_i^{(1)} \ot  {\partial r^{23}
\over \partial q_i}(q, w-cv),   \  \xi (z) \ot \eta (w )
\ot \zeta (cv) \ra  dz dw dv\\
&=&\la \zeta , \ \oint_{C} \sum_{i} \xi_{i}(z) {\partial \over \partial q_i}
(R_{q}+\half I )\eta dz \ra \\
&=& \la \zeta ,  \calx_{j^*\xi} (R_{q}\eta ) \ra ,
\ee
and similarly
\be
&&(\frac{1}{2\pi i})^{3} \oint_{C}  \lim_{\ct 1-} \oint_{C}
\oint_{C}  \la \sum _i \, h_i^{(2)} \ot  {\partial r^{31}
\over \partial q_i}(q, cv -z),   \  \xi (z) \ot \eta (w )
\ot \zeta (cv) \ra  dz dw dv\\
&=&- \la \zeta  ,  \calx_{j^*\eta } (R_{q}\xi  ) \ra .
\ee
Lastly, it follows from Equation (\ref{eq:3.15c}) that
\be
&&(\frac{1}{2\pi i})^{3} \oint_{C}  \lim_{\ct 1-} \oint_{C}
p.v. \oint_{C} \la  \sum_i \, h_i^{(3)} \ot  {\partial r^{12}
\over \partial q_i}(q, z -w),   \  \xi (z) \ot \eta (w )
\ot \zeta (cv) \ra  dz dw dv\\
&=& \frac{1}{2\pi i}
 \oint_{C} \la \zeta (v),      \sum_{i} {\partial \over \partial q_i}
\la R_{q}\xi , \eta   \ra  h_{i} \ra   dv\\
&=& \la \zeta ,  \      d  \la R_{q}\xi , \eta   \ra \ra .
\ee
Assembling the calculation, using the fact that  $r$ satisfies
(CDYBE), we conclude that $R_q$ satisfies (mCDYBE).
The assertion that $R_q \in \call (L\frakg^{*}  , L\frakg)^{H}$ 
now follows from Equation (\ref{eq:4.11}) and Lemma \ref{lem:3.13a}(1). \qed

According to Proposition \ref{pro:R}
and  Corollary \ref{cor:R},  we can use $R$ to equip
  $T^*\frakh^{*} \times L\frakg^* $
 with a Lie algebroid structure, and therefore 
$T\frakh^{*} \times L\frakg$ admits  the   Lie-Poisson structure.
On the other hand, consider   $T^*\frakh^{*}$ with 
the  canonical cotangent  symplectic structure, $\frakg^*$
  with the  plus Lie Poisson structure,
 and equip $T^*\frakh^{*} \times \frakg^*$
with the product Poisson structure. According to Example 2.1,
 this product structure is just the Lie-Poisson structure on the dual
vector bundle $T^*\frakh^{*} \times \frakg^* $, when $T\frakh^{*} \times \frakg$
is the product  Lie algebroid.
 In the next proposition, we are going
to establish a Poisson map from $T^*\frakh^{*} \times \frakg^* $ to
$T\frakh^{*} \times L\frakg$. This  essentially enables us
to describe certain finite-dimensional symplectic
leaves of $T\frakh^{*} \times L\frakg$, which are simply the image
of $T^*\frakh^{*} \times \calo$  under this map
for    coadjoint orbits $\calo\subset \frakg^*$.
In order to do so, we need an equation somewhat intermediate
between (CDYBE) and (mCDYBE) which involves both 
$(r_{-}^{\#}(q))^{*}$ and $R_q$:
\begin{eqnarray}
\label{eq:4.17}
&&[(r_{-}^{\#}(q))^{*} \xi, \ (r_{-}^{\#}(q))^{*} \eta ]
-(r_{-}^{\#}(q))^{*} I^{-1} ([R_{q}\xi ,  I\eta ]+[I\xi , R_{q}\eta ]) \nonumber\\
&& +\calx_{j^*\xi } ((r_{-}^{\#}(q))^{*} \eta )-
\calx_{j^*\eta  } ((r_{-}^{\#}(q))^{*} \xi  ) 
+d[\la R_{q}\xi , \eta \ra ]=0, \ \ \ \ \xi, \eta \in L\frakg^* .
\end{eqnarray}
The derivation of this equation makes use of Equations
(\ref{eq:3.14a})-(\ref{eq:3.14c})  and (\ref{eq:3.15a})-(\ref{eq:3.15c})
 with $\zeta (v) $
replaced by $\zeta \in \frakg^*$ and with $v=0$. As the calculation
is similar to the proof of Proposition \ref{pro:Rmdybe}, we shall omit
the details.

\begin{thm}
\label{pro:Poisson}
The map $\rho : T^*\frakh^{*} \times \frakg^{*} \lon
 T\frakh^{*} \times L\frakg$ 
given by
\begin{equation}
\label{eq:4.190}
(q, p, \xi )\lon (q, -i^*\xi, p+r_{-}^{\#}(q)\xi), \ \ q\in \frakh^{*},
p \in \frakh, \xi \in \frakg^* , 
\end{equation}
is an $H$-equivariant Poisson map, where $H$ acts on
$T^*\frakh^{*} \times \frakg^{*} $ and 
 $T\frakh^{*} \times L\frakg$ by acting on the second factors
by coadjoint and adjoint actions respectively,
$i: \frakh \lon \frakg$ is the natural
inclusion, and $i^* :\frakg^*\lon \frakh^*$ is the dual map. 

In other words, $\rho$ is an $H$-equivariant  realization  in the  dynamical
Lie algebroid $T^* \frakh^{*} \times L\frakg^*$ in the sense
of Definition \ref{def:3.1}.
\end{thm}
\pf In order to show that $\rho$ is  a Poisson map, it is
enough to check that the dual map $\rho^* :T^*\frakh^{*} \times L\frakg^*
\lon T\frakh^{*} \times \frakg$ is a  morphism of Lie algebroids.
By direct calculation, we  have 
$\rho^* (q, p, \xi )=(q, j^* \xi , -p +(r_{-}^{\#}(q))^* \xi )$,
$q\in \frakh^{*}, \ p\in \frakh$, $\xi \in L\frakg^*$. There are
two conditions to check. First, we have to show that 
$a\smalcirc  \rho^* =a_{*}$, where $a: T\frakh^{*}  \times \frakg
\lon T\frakh^*$ is anchor map of the trivial  Lie algebroid.
From the definition of the various quantities, this is trivial.
Secondly, we have to check that the induced map
on sections preserve the Lie algebroid brackets.
To do so, it is enough to verify that this
is the case for brackets between constant sections. Thus
we have to check that
\begin{enumerate}
\item $\rho^* [(h_{1}, 0), \ ( h_{2}, 0)]=[\rho^* (h_{1}, 0), \ \rho^* ( h_{2}, 0)], \ \ \forall h_{1}, h_{2}\in \frakh$;
\item $\rho^*  [(h , 0),  \ ( 0, \xi )]= [ \rho^* (h , 0),  \ \rho^*  ( 0, \xi )], \ \ \forall h\in \frakh, \xi \in L\frakg^*$;
\item $\rho^*  [(0, \xi ),  \ ( 0, \eta )]= [ \rho^*  (0, \xi ),  \ \rho^* 
 ( 0, \eta  ) ], \ \ \forall \xi , \ \eta  \in L\frakg ^* $;
\end{enumerate}

For (1), the equality follows because $\frakh$ is Abelian.
For (2), we have 
\be
&&\rho^* [(h , 0),  \ ( 0, \xi )]\\
&=& (-j^* ad_{h}^{*}\xi , \  -(r_{-}^{\#}(q))^* ad_{h}^{*} \xi )\\
&=&(0, \  -(r_{-}^{\#}(q))^* ad_{h}^{*}\xi ).
\ee
 as $j^* ad_{h}^{*} \xi=0$.
On the other hand, 
\be
&&[ \rho^* (h , 0),  \ \rho^*  ( 0, \xi ) ]\\
&=&[(0, -h), (j^* \xi , (r_{-}^{\#}(q))^* \xi )]\\
&=&(0, -ad_{h}(r_{-}^{\#}(q))^* \xi ).
\ee
 Hence the result follows
from Lemma \ref{lem:3.13a}. For (3), we have
\be
&&\rho^*  [(0, \xi ),  \ ( 0, \eta )]\\
&=&\la j^* I^{-1} ([R\xi , I \eta ]+[I\xi , R \eta ]), \ -d\la R\xi , \eta \ra 
+ (r_{-}^{\#}(q))^* I^{-1} ([R\xi ,  I\eta ]+[I\xi , R \eta ] ) \ra\\
&=&( 0, -d\la R\xi , \eta \ra +(r_{-}^{\#}(q))^* I^{-1}
 ([R\xi ,  I\eta ]+[ I\xi , R \eta ]))   \ \ \ \ \mbox{(by Lemma \ref{lem:3.13a})}.
\ee

On the other hand,

\be
&& [\rho^*  (0, \xi ),  \ \rho^*  ( 0, \eta )]\\
&=&[(j^* \xi , (r_{-}^{\#}(q))^* \xi ), \ (j^* \eta , (r_{-}^{\#}(q))^* \eta )]\\
&=&([j^* \xi , j^* \eta ], \  \calx_{j^* \xi } (r_{-}^{\#}(q))^* \eta
-\calx_{j^* \eta }(r_{-}^{\#}(q))^* \xi 
+[ (r_{-}^{\#}(q))^* \xi , \  (r_{-}^{\#}(q))^* \eta ]).
\ee
Therefore the equality $\rho^* [(0, \xi ),  ( 0, \eta )]
=[\rho^* (0, \xi ),\ \rho^* ( 0, \eta )]$
follows from the commutativity of $\frakh$ and Equation (\ref{eq:4.17}). \qed

Following the notations in Section 3 (Equations (\ref{eq:3.3}-\ref{eq:3.5})),
we  have

\begin{equation}
\label{eq:4.3a}
        L : T^*\frakh^{*} \times \frakg^*  \to L\fg, \ \ \ \ 
L(q, p, \xi )= Pr_{2}\smalcirc \rho (q, p, \xi )=  p+r_{-}^{\#}(q)\xi;
\end{equation}
\begin{equation}
\label{eq:4.4a}
        \tau :T^*\frakh^{*} \times \frakg^* \to T \frakh^* ,  \ \ \ \ 
\tau (q, p, \xi )= Pr_{1}\smalcirc \rho (q, p, \xi )=  (q, -i^{*}\xi );
\end{equation}
and
\begin{equation}
\label{eq:4.5a}
        m :T^*\frakh^{*} \times \frakg^* \to   \frakh^* ,  \ \ \ \ 
 m (q, p, \xi )=  p \smalcirc \tau (q, p, \xi )= q. 
\end{equation}

\begin{defi}
A function on $L\fg$ is said to be smooth on  $L\fg$ 
if for each $X\in L\fg$, the derivative
$df (X)\in L\fg^*$ (recall that $df (X)$ is defined
as a linear functional on  $L\fg$ through the
relation  $\left.\frac{d}{dt}\right|_{t = 0} f(X + tY) =
 df(X) (Y ), \ \forall X, Y\in L\frakg$).
\end{defi}

Combining Theorem  \ref{pro:Poisson} with
Propositions \ref{pro:3.6}, \ref{pro:3.8},
 we are lead to the   following 

\begin{thm}
 Assume that $r$ is a     classical
dynamical $r$-matrix with  spectral parameter. Then
\begin{enumerate}
\item  $L:  T^*\frakh^{*} \times \frakg^*\lon  L\frakg$, 
$(q, p, \xi )\lon  p+r_{-}^{\#}(q)\xi $
satisfies
\begin{eqnarray}
        & & \{ L^*f, \ L^*g \}(x)   \label{eq:LLL} \\
        & = & \langle L(x), \ -ad_{R_{m(x)}(df(L(x)))}^{*} dg(L(x)) +
                ad_{R_{m(x)}(dg(L(x)))}^{*}  df(L(x))\rangle  \nonumber \\
        & & + \langle (\calx_{\tau(x)}R)(df(L(x))), \ dg(L(x))\rangle,\quad\quad
                \mbox{for }  x=(q, p, \xi ) \in T^*\frakh^{*} \times \frakg^*,  \nonumber
\end{eqnarray}
and  all smooth functions $f, \ g$ on $L \frakg$.
\item
If  $\HH =L^* f, \ f\in I(L\frakg)$,  then under  the  
flow  $\phi_{t}$ generated by the Hamiltonian $\HH$, we have the
following   quasi-Lax type equation:
\begin{equation}
\frac{dL(\phi_{t})}{dt}=
[R_{m(\phi_{t})} (df(L (\phi_{t}))), L(\phi_{t})] 
-(\calx_{\tau ( \phi_{t} )}  R )(df(L(\phi_{t}))).
\end{equation}
\end{enumerate}
\end{thm}

\begin{rmk}
\label{rmk:4.9}
{\em In the first part of the above theorem, we have
restricted ourselves to smooth functions on $L\fg$
with derivatives in the restricted dual $L\frakg^*$.
However, we can easily extend the calculation to include
linear functions of the form $l_{\xi}(X)=\la \xi , X(z)\ra$,
where $\xi \in \frakg^*$ and $X\in L\frakg$.
For these functions, the derivative  $dl_{\xi}(X)=
\delta (z- \cdot )\xi$ is in the singular part of
$(L\frakg )^*$, where $\delta$ is the delta
function. In particular, we obtain the  St. Petersburg
type formula:
\begin{eqnarray}
\{L (z) \  \stackrel{\ot}{,}\  L (w)\}&=&-[r^{12}(q, z-w), L^{1}(z)]
+[r^{21}(q, w-z), L^{2}(w)] -\calx_{i^*\xi}r(q, z-w)\\
&=&-[r^{12}(q, z-w),  L^{1}(z)+ L^{2}(w) ]-\calx_{i^*\xi}r(q, z-w)
\end{eqnarray}
by calculating with such linear functions. Here,
$L(z):  T^*\frakh^{*} \times \frakg^*\lon  L\frakg$
is defined by $L(z)(q, p, \xi )=L(q, p, \xi )(z)$ and it is
understood that $L^{1} (z)=L(z)\ot 1$ and 
$L^{2}(w)=1\ot L(w)$ in the above
formula are evaluated at $(q, p, \xi )$.
}
\end{rmk}

In the rest of the section, we shall   consider the case where
$\fg$ is a simple Lie algebra over $\complex$ with Killing form
$(\cdot , \cdot )$ and we shall take $\frakh$ to be a 
fixed Cartan subalgebra.

Let $Q$ be the quadratic function
\begin{equation}
Q(X)=\half \oint_C (X(z), X (z))\frac{dz}{2\pi i z}, \ \ \ \forall X\in L\frakg,
\end{equation}
where $C$ is a small circle around the origin.
Clearly, $Q$ is an ad-invariant function on $L\frakg$.

\begin{defi}
\label{defi:spin}
Assume that $r$ is a     classical
dynamical $r$-matrix with  spectral parameter.
The Hamiltonian system on $ T^*\frakh^{*} \times \frakg^*$  generated
by the Hamiltonian function:
\begin{equation}
\HH  (q, p, \xi)= (L^*Q)  (q, p, \xi)= \half \oint_C (L(q, p, \xi), \  L(q, p, \xi))
\frac{dz}{2\pi iz}
\end{equation}
is called the  spin Calogero-Moser system associated to the
dynamical r-matrix $r$.
\end{defi}

In  \cite{EV}, Etingof and Varchenko 
obtained a complete classification of classical
dynamical $r$-matrices  (which satisfy
Equations (\ref{eq:4.1}-\ref{eq:cdybe}))
for simple Lie algebras.
 Up to  gauge transformations, they obtained
canonical forms of the three types  (rational, trigonometric
and elliptic) of dynamical $r$-matrices.
 For  each of these dynamical r-matrices,
one can associate a spin Calogero-Moser system on $ T^*\frakh^{*} \times \frakg^*$.  We will  list all of them below (see Remark \ref{rmk:4.2}(1)).
First, let us fix some notations. Let $\frakg =\frakh \oplus 
\sum_{\al \in \Delta}\frakg_{\al}$ be the root
space decomposition. For any positive
root $\alpha \in \Delta_{+}$, fix basis  $e_{\al}\in \frakg_{\al}$
and $e_{-\al}\in \frakg_{-\al}$  which are dual with respect to
$(\cdot , \cdot )$.
Fix also  an orthonormal  basis $\{h_{1}, \cdots , h_{N}\} $
of $\frakh$, and write $p=\sum_{i=1}^N p_{i}h_{i}$,
$\xi_{i}=\la \xi , h_{i}\ra$,  and $\xa= \la \xi , e_{-\al} \ra $,
 for $p\in \frakh$ and $\xi \in \frakg^*$. Then
$I\xi =\sum_{i=1}^N \xi_{i} h_{i}+\sum_{\al \in \Delta} \xa e_{\al} \in \frakg$.\\\\

{\bf I. {{ Rational case}}}

\be
&&r(q, z) \, = \,  {\Omega  \over z} \,
+ \sum _{\alpha  \in \Delta'}  \,
{1\over (\al, q) } \,
e_\al \otimes  e_{-\al},\\
&& \HH (q, p, \xi )=\half \sum_{i=1}^{N}p_{i}^{2}-\half \sum_{\alpha \in \Delta'}
\frac{\xa \ax}{(\alpha , q)^2},\\
&&L(q, p, \xi )(z)=p+\frac{I\xi}{z}+\sum_{\alpha \in \Delta'}
\frac{\xa}{(\alpha , q)}e_\al ,
\ee
where  $\Delta' \subset \Delta$ is a set of roots
 closed with respect to the addition and multiplication by $-1$.\\\\

{\bf II.  {{ Trigonometric  case}}}

\be
&&r(q, z)\,=\, (\mbox{cot}\, z +\third z ) \,\sum_{i=1}^N h_i\otimes h_i\,
+ \, \sum _{\al \in \Delta (\Pi' )} \,{ \mbox{sin}\, ((\al, q) + z)
\over
 \mbox{sin}\, (\al, q) \,  \mbox{sin}\, z} 
e^{\third z (\alpha , q)}\,e_\al \otimes  e_{-\al} \\
&&\ \ \ \ \ \ \ \ \ \ + \, \sum _{\al \in \Delta_{+}-\Delta (\Pi'  )}\,{ e^{-iz}
\over
 \mbox{sin} \,z} e^{\third z (\alpha , q)} \,e_\al \otimes  e_{-\al} 
+ \, \sum _{\al \in \Delta_{-}-\Delta (\Pi'  )}\,{ e^{iz}
\over
 \mbox{sin}\, z}e^{\third z (\alpha , q)}\,e_\al \otimes  e_{-\al}, \\
&& \HH (q, p, \xi )=  \half \sum_{i=1}^{N}p_{i}^{2}-\half
\sum _{\al \in \Delta (\Pi'  )}(\frac{1}{\sin^{2}(\al, q)}-\frac{1}{3})\xa \ax
-\frac{5}{6}\sum _{\al \in \Delta-\Delta (\Pi'  )}\xa \ax ,\\
&&L(q, p, \xi)(z)=p+(\cot{z}+\frac{1}{3} z) \sum_{i=1}^N \xi_{i} h_{i} +
\sum _{\al \in \Delta (\Pi'  )}{ \mbox{sin}\, ((\al, q) + z)
\over
 \mbox{sin}\, (\al, q) \,  \mbox{sin}\, z}e^{\frac{1}{3}z (\al , q)}\xa e_\al \\
&&\ \ \ \  \ \ \ \ \ \  \ \ + \sum _{\al \in \Delta_{+}-\Delta (\Pi'  )}\,{ e^{-iz}
\over
 \sin{ \,z}} e^{\frac{1}{3}z (\al , q)}\xa e_\al
+\sum _{\al \in \Delta_{-}-\Delta (\Pi'  )}\,{ e^{iz}
\over
 \mbox{sin}\, z}e^{\frac{1}{3}z (\al , q)}\xa e_\al .
\ee
Here $\Delta =\Delta_{+}\cup \Delta_{-}$
is a polarization of $\Delta$, 
 $\Pi' $ is a  subset of the set of simple  roots,
and  $\Delta (\Pi'  )$
denotes the set of all roots which are linear
combinations of roots from $\Pi' $.\\\\

{\bf  III.  {{  Elliptic case }}}

\be
&&r(q, z)\,=\, \zeta (z)\sum_{i=1}^N h_i\otimes h_i\,
- \,\sum_{\al \in \Delta} l(q, z ) e_\al \otimes e_{-\al}, \\
&&\HH (q, p, \xi)=  \half \sum_{i=1}^{N}p_{i}^{2}-\half
\sum _{\al \in \Delta }\calp ((\al , q ))\xa \ax ,\\
&&L(q, p, \xi)(z)=p+\zeta (z)
\sum_{i=1}^N \xi_{i} h_{i} -\sum _{\al \in \Delta } l((\al , q ), z) \xa e_{\al},
\ee
where $\zeta (z)=\frac{\sigma' (z)}{\sigma (z)}$, $\calp (z)
=-\zeta' (z)$,
 $l(w, z)=-\frac{\sigma (w+z)}{\sigma (w)\sigma (z)}$, and
$\sigma (z)$ is the Weierstrass $\sigma$ function of periods 
$2\omega_{1}, 2\omega_2$.\\\\\\ 

\begin{rmk}
\label{rmk:4.2}
{\em \begin{enumerate}
\item In the trigonometric case and the elliptic case,
the classical   dynamical $r$-matrices with spectral 
parameter which we used above are gauge equivalent to  those
 in \cite{EV}. If we had used the canonical
forms given in \cite{EV}, the Hamiltonians
of the associated spin systems will have additional
terms which depend on $i^* \xi$. The same remark
also applies to the most general dynamical
r-matrix which one can obtain by using gauge
transformations. However, as will be
evident in the next result, these
additional terms do not give rise to any
new systems upon reduction.
\item In the rational  and  trigonometric case above, the spin systems that
we have here are in one-to-one correspondence with some subsets
of the root system. Thus we have as many spin systems as these
special subsets.
\item The reader should note that the  ${\frak so} (N)$ models in \cite{BAB2}
are   different from ours.
\end{enumerate}
}
\end{rmk}

We conclude this section with
the following result which prepares the way for the
construction of associated integrable models in the next section.

\begin{thm}
\label{thm:aaa}
 The Hamiltonians of the spin Calogero-Moser
systems are invariant under the canonical $H$-action
on $ T^*\frakh^{*} \times \frakg^*$:
\begin{equation}
\label{eq:action}
x\cdot (q, p, \xi )=(q, p, Ad_{x^{-1}}^* \xi ), \ \ \forall x \in H, \ 
(q, p, \xi )\in T^*\frakh^{*} \times \frakg^* 
\end{equation}
 with momentum map $J: T^*\frakh^{*} \times \frakg^* \lon \frakh^{*}$
given by 
\begin{equation}
\label{eq:J}
J(q, p, \xi ) =i^*\xi.
\end{equation}
 If $\Sigma$ denotes
the set defined by $\calx_{i^*\xi}R=0$, then
$\Sigma=J^{-1}(0)$ in the trigonometric and elliptic
cases, while $\Sigma =J^{-1} ((\Delta ')^{\perp})$ in
the rational case. Thus in each case, $\Sigma$ is invariant
under the dynamics and we have $\frac{dL}{dt}=[R (M) , L]$
on $\Sigma$, where $M(q, p, \xi )(z) =L (q, p, \xi )(z)/z$. 
\end{thm}

\begin{rmk}
{\em Note that  $J^{-1}(0)$ is not a Poisson submanifold of 
$T^*\frakh^{*} \times \frakg^* $, otherwise
the corresponding subsystem on $J^{-1}(0)$  would  have a
  natural collection of Poisson commuting integrals and there
would have no need to use reduction to 
construct the associated integrable flows. 
}
\end{rmk}

\section{Integrable Spin \CM systems}

In this section, we shall carry out the reduction procedure
outlined in Section 3 to the spin \CM systems. As a result,
we obtain a  new family of integrable systems, which we
call {\em integrable spin \CM systems}. For 
$\frakg =\frak{sl}(n, \complex )$, the usual \CM systems
as well as their spin generalizations (in the sense of
Gibbons and Hermsen \cite{GH}) appear as subsystems
of what we have on special symplectic leaves
of the reduced Poisson manifold. However, for other simple Lie algebras,
 the usual \CM systems without spin cannot be realized in this 
fashion, as we shall explain below.
 
Our first task below is to construct an $H$-equivariant 
map $g$ which allows us to construct the equations
of motion in Lax pair form for the reduced Hamiltonian
$\HH_0$. 

For any root $\alpha \in \Delta$, recall that the coroot
$h_{\alpha}$ is the element in $\frakh$ corresponding
to $2\frac{\alpha}{(\alpha , \alpha )}$ under
the isomorphism between $\frakh$ and $\frakh^*$ induced by the
Killing form $(\cdot , \cdot )$.
I.e., for any $\beta \in \frakh^*$, $\beta (h_{\alpha})=
2\frac{(\beta , \alpha )}{(\alpha, \alpha )}$. Therefore,
if we fix a simple system $\Pi =\{\alpha_{1}, \cdots , \alpha_{N}\}
\subset \Delta$, we have a  basis of $\frakh$ given by
the fundamental coroots $h_{\alpha_{1}}, \cdots , h_{\alpha_{N}}$.
In particular, the entries of the Cartan matrix  $A=(A_{ij} )$
is given by $A_{ij}=\alpha_{j}(h_{\alpha_{i}})$. Let
$\omega_{1}, \cdots , \omega_{N}$ be the fundamental
weights, i.e., the dual basis of $h_{\alpha_{1}}, \cdots , h_{\alpha_{N}}$ 
in $\frakh^*$.
Then it is clear that
\begin{equation}
\label{eq:5.1}
 \alpha_{i}=\sum_{j=1}^{N} A_{ji}\omega_{j}.
\end{equation}

We shall denote by $C=(C_{ij})$ the inverse of the Cartan matrix. 
Clearly, we have $C_{ij}\in  {\Bbb Q}, \ \forall i, j$.
Consider  the  open  submanifold  of $\frakg^*$:
\begin{equation}
\calu =\{\xi\in \frakg^* |\xi_{\alpha_i} = \la \xi , e_{-\alpha_i}\ra\neq 0, \ i=1, \cdots , N\}.
\end{equation}
It is clear that $\calu$ is stable under the coadjoint action
of $H$ (considered as a subgroup of $G$).
Our next aim is to construct a map $g: \calu \lon H$ with the
property that
\begin{equation}
\label{eq:g}
 g(Ad_{h^{-1}}^{*} \xi )=h \cdot  g(\xi ), \ \ \ \ \ \forall  h\in H.
\end{equation}
In other words, $g$ is equivariant, where $H$ acts on itself  by
left translation.

For the sake of convenience, below we will identify $\frakg^*$ with
 $\frakg$ by the Killing form and  identify $\calu$ with
the open submanifold $\{ \xi=\sum_{i=1}^N \xi_{i} h_{i}+
\sum_{\al \in \Delta} \xa e_{\al}  \in \frakg|  
\xi_{\alpha_i}\neq 0, \ i=1, \cdots , N\}$ of $\frakg$.
Thus the coadjoint action becomes the  adjoint action  and Equation 
(\ref{eq:g}) becomes
\begin{equation}
\label{eq:g1}
 g(Ad_{h} \xi )=h \cdot  g(\xi ), \ \ \ \ \ \forall  h\in H.
\end{equation}

Since $H$ is generated by a small neighborhood of $1$, it is
 sufficient for $g$ to satisfy Equation (\ref{eq:g1})
for  $h \in U \subset H$, where $U$ is sufficiently
small so that  the map $\log : U\lon \frakh$ inverse
to the exponential map is well defined. 
Indeed, for all $h\in U$,  we have
\begin{eqnarray}
&&\log{h}  =\omega_{1}(\log{h} )h_{\alpha_1}+\cdots +
\omega_{N}(\log{h} )h_{\alpha_N},\ \ \mbox{ and}\label{eq:5.3a} \\
&&Ad_{h} e_{\alpha }=\chi_{\alpha} (h) e_{\alpha }, \ \ \ \
\chi_{\alpha} (h)=e^{\alpha (\log{h} )}. \label{eq:5.3b}
\end{eqnarray}

Note that if $g_{i}: \calu \lon H, \ i=1, \cdots , N$,
 satisfies that
\begin{equation}
\label{eq:gi}
g_{i} (Ad_{h} \xi )= \exp{ (\omega_{i}( \log{h})h_{\alpha_{i}}} )g_{i}  (\xi ),
\ \ \forall h\in U,\   \xi \in \calu,  \ \ \ i=1, \cdots , N, 
\end{equation}
 then $g=g_{1}\cdots g_{N}$ will have the desired property in
Equation (\ref{eq:g1}).
Next we shall  seek $g_{i}$ in the form
\begin{equation}
\label{eq:5.5}
g_{i}(\xi )=\exp{(\phi_{i} (\xi )h_{\alpha_{i}} )} ,
\end{equation}
where $\phi_{i}$ is a function on $\calu$.
In order for $g_{i}$ to satisfy Equation (\ref{eq:gi}), it is enough
that
\begin{equation}
\label{eq:phii}
\phi_{i} (Ad_{h} \xi )=\phi_{i} ( \xi )+\omega_{i} (\log{h} ).
\end{equation}
Let $\psi_{i}(\xi )=  e^{\phi_{i}(\xi )}$. Then Equation (\ref{eq:phii})
translates into
\begin{equation}
\label{eq:psii}
\psi_{i} (Ad_{h} \xi ) =\chi_{i}(h)\psi_{i}(\xi ), \  \ \ 
\chi_{i}(h)= e^{\omega_{i}(\log{h} )}.
\end{equation}
That is, $\psi_{i}$ is  a semi-invariant with character $\chi_{i}$.
In what follows,  we shall fix a branch of the logarithmic
function. We shall seek  $\psi_{i}$ of  the form 
\begin{equation}
\label{eq:5.8}
\psi_{i} (\xi )=\prod_{j=1}^{N} \xi_{\alpha_{j}}^{n_{ij}}, \ \ \ \forall
 \xi\in \calu. 
\end{equation}
Then by Equations (\ref{eq:5.3a}-\ref{eq:5.3b}), 
\be
&&\psi_{i} (Ad_{h} \xi )\\
&=&\prod_{j=1}^{N} (\chi_{\alpha_{j} }(h)\xi_{\alpha_{j} })^{n_{ij}}\\
&=&(\prod_{j=1}^{N}\chi_{\alpha_{j} }(h)^{n_{ij}} )\psi_{i}(\xi )\\
&=&(\prod_{j=1}^{N} e^{n_{ij}\alpha_{j} (\log{h} )})\psi_{i}(\xi ),
\ \ \ h\in U, \ \xi \in \calu.
\ee
Therefore, in order to satisfy Equation (\ref{eq:psii}), it suffices
to pick $n_{ij}$ so that 
$\omega_{i}=\sum_{j=1}^{N} n_{ij}\alpha_{j}$.
But from the relation in Equation (\ref{eq:5.1}), we must have  
 $n_{ij}=C_{ji}$, i.e.,
\begin{equation}
\label{eq:5.9}
\psi_{i} (\xi )= \prod_{j=1}^{N} \xi_{\alpha_{j}}^{C_{ji}}
\end{equation}
and
\begin{equation}
g_{i}(\xi )=\exp{ (\sum_{j=1}^{N} C_{ji} \log{\xi_{\alpha_{j}}})
h_{\alpha_{i}}}.
\end{equation}

Consequently, we have

\begin{thm}
\label{thm:g}
The formula
\begin{equation}
\label{eq:5.11}
g(\xi ) =\exp{(\sum_{i=1}^{N}
\sum_{j=1}^{N}(C_{ji} \log{\xi_{\alpha_{j}}})h_{\alpha_{i}})}
\end{equation}
defines an $H$-equivariant map $g: \calu\lon H$.
\end{thm}
 
Consider the Poisson submanifold $T^* \frakh^* \times \calu $
of  $T^* \frakh^* \times \frakg^*$. Clearly,   the 
$H$-action defined by Equation  (\ref{eq:action}) 
induces a Hamiltonian action on $T^* \frakh^* \times \calu $
and therefore the moment map $J: \  T^{*}   \frakh^{*} \times \calu \lon \frakh^* $
is given by restriction of the one in Equation (\ref{eq:J}).
 Hence $J^{-1}(0)=T^* \frakh^* \times (\hp \cap \calu )$,
and therefore   we have
\begin{equation}
\label{eq:574}
\calx_{v}R=0, \ \ \ \ \ \forall v\in \tau (J^{-1}(0) ).
\end{equation}

 Thus according to  Theorem \ref{pro:Poisson}, Theorem  \ref{thm:aaa}, 
Theorem \ref{thm:g}  and
 Equation (\ref{eq:574}), 
  we conclude that Assumptions A1--A4 in Section 3 are all
satisfied and therefore we can now apply the reduction procedure
of  Section 3 to our situation.
We first characterize the reduced
space using  the following:

\begin{thm}
\label{thm:5.12}
The quotient space $J^{-1}(0)/H\cong T^* \frakh^* \times (\hp \cap \calu )/H$
is analytic and can be identified with $T^* \frakh^* \times \frakg^*_{red}$,
where $\frakg^*_{red}$ is the affine subspace $\epsilon +\sum_{\alpha
\in \Delta -\Pi} \complex e^*_{\alpha }$ and $\epsilon =\sum_{i=1}^{N}
e^*_{\alpha_j}$, where $\{e^*_{\alpha }|\alpha \in \Delta\}$ denotes
the dual  vectors in $\frakg^*$ corresponding to
$\{e_{\alpha }|\alpha \in \Delta\}$ in $\frakg$.
\end{thm}
\pf  It is simple to see that the action of $H$  on $\hp \cap \calu $ is locally
free.
Moreover, each $H$-orbit through $\hp \cap \calu $ has  exactly one
intersection with $\frakg^*_{red}$. To see this, recall 
that the simple system has the following property, namely,
if $\alpha \in \Delta $, there exist integers $m_{\alpha}^{i}\ 
(1\leq i\leq N)$ either all nonnegative or all nonpositive,
such that $\alpha =\sum_{i=1}^{N}m_{\alpha}^{i} \alpha_{i}$.
Hence for a given $\xi \in \hp \cap \calu$, if we let $h=
g(\xi )^{-1}$, then
$Ad_{h^{-1} }^* \xi =\sum_{\alpha \in \Delta}\xi_{\alpha} \chi_{\alpha} (h)
e^*_{\alpha }=\epsilon +\sum_{\alpha \in \Delta -\Pi}
(\xi_{\alpha} \prod_{i=1}^{N} \xi_{\alpha_i}^{-m_{\alpha}^{i}})
e^*_{\alpha }\in  \frakg^*_{red}$.
Hence we can identify $(\hp \cap \calu )/H$ with
$\frakg^*_{red}$. \qed

\begin{rmk}
\label{rmk:5.1}
{\em Indeed $s_{\alpha} =\xi_{\alpha} \prod_{i=1}^{N} \xi_{\alpha_i}^{-m_{\alpha}^{i}}$, $\alpha \in \Delta -\Pi$, are a set of $H$-invariant 
functions on $\hp \cap \calu $, which can be used as a coordinate
system for $\frakg^*_{red}$. 
If $s \in \frakg^*_{red}$, we then may   write $s =\sum_{\alpha \in
\Delta} s_{\alpha} e^*_{\alpha}$ with $s_{\alpha_{i}}=1, \
i=1, \cdots , N$.}
\end{rmk}

By Poisson reduction \cite{MR}, the reduced manifold
$T^* \frakh^* \times \frakg^*_{red}$ has a unique Poisson
structure which is a  product structure, where the
second factor $\frakg^*_{red}$ is being equipped with  the reduction
(at $0$) of  the Lie-Poisson  structure on $\calu$ by the
$H$-coadjoint action.  The Poisson brackets between
the coordinate functions $s_{\alpha}$ on $\frakg^*_{red}$
can be obtained by a straightforward but tedious
 computation. We shall leave
the details to interested reader.
  Now, the symplectic leaves of $\frakg^*_{red}$ 
are the symplectic reduction  of $ \calo \cap \calu$ at $0$,
where $\calo \subset \frakg^*$ is a coadjoint orbit \cite{MR}.
In other words, any symplectic leaf of $\frakg^*_{red}$
is of the form $(\calo \cap \calu \cap \hp ) / H$, and we
shall  
denote this by  $\calo_{red}$. Obviously, $\calo_{red}$ is a
symplectic manifold of dimension $\mbox{dim} \calo- 2 N$, where
$N$ is the rank of the Lie algebra $\frakg$.
Consequently, the  symplectic leaves of $T^* \frakh^* \times \frakg^*_{red}$
are of the form $T^* \frakh^* \times \calo_{red}$, which is of
dimension equal to $\mbox{dim} \calo$.

Accordingly, if $\HH$ is the Hamiltonian of one of the spin \CM systems
in Section 4, and $L$ is the   corresponding Lax operator, there
exists uniquely determined Hamiltonian function $\HH_0$ and Lax operator
$L_0$ on the reduced Poisson manifold $T^* \frakh^* \times \frakg^*_{red}$
such that $\ppi^*   \HH_0 =\HH |_{T^* \frakh^* \times (\hp \cap \calu )}$
and $L_{0}\smalcirc \ppi =\tilde{L}|_{T^* \frakh^* \times (\hp \cap \calu )}$.
Here,  $\ppi : T^* \frakh^* \times (\hp \cap \calu )\lon T^* \frakh^* \times
\frakg^*_{red}$ is the natural projection given by
\begin{equation}
\label{eq:5.13}
\ppi  (q, p ,\xi )=(q, p,  \ \epsilon +\sum_{\alpha \in \Delta -\Pi}
(\xi_{\alpha} \prod_{i=1}^{N} \xi_{\alpha_i}^{-m_{\alpha}^{i}})
e^*_{\alpha }),
\end{equation}
and $\tilde{L} :T^* \frakh^* \times \calu \lon L\frakg$ is given
by 
\begin{equation}
\tilde{L} (q, p, \xi )= Ad_{g(\xi )^{-1}} L(q, p, \xi ).
\end{equation}

We can now state the main result of the paper.

\begin{thm}
\label{thm:5.13}
Let  $\HH$ be  the Hamiltonian of a spin \cm system with  Lax
operator ${L}$.
And let  $ \tilde{R}:  T^* \frakh^* \times  \calu
\lon \call (L\frakg^* , L\frakg )$ be  the map
as defined by Equation (\ref{eq:3.17}), which is obtained from
 $R$ by applying the gauge transform given by $g$ in Theorem \ref{thm:g},
and $R_{0}$ the induced map on $ T^* \frakh^* \times \frakg^*_{red}$
in the sense that $R_{0}\smalcirc \pi =\tilde{R}|_{T^* \frakh^* \times (\hp \cap \calu )}$.
  Then the Hamiltonian   system generated by the induced function $\HH_0$ on the
 reduced Poisson manifold 
$T^* \frakh^* \times  \frakg^*_{red}$ admits a
Lax operator $L_{0}:  T^* \frakh^* \times  \frakg^*_{red} \lon 
L\frakg$ satisfying the following properties:
\begin{enumerate}
\item  For any  smooth functions  $ f_1, \ f_2 $
on $ L\fg$,
\begin{eqnarray}
        & & \{L_{0}^* f_1, \  L_{0}^* f_2\} (\xx) \nonumber \\
        & = & - \langle L_{0}(\xx), \ ad^{*}_{R_{0}^{*}(\xx)
                (df_2(L_{0}(\xx)))} df_1(L_{0}(\xx))  \nonumber \\
        & & \qquad + ad_{R_{0}(\xx)(df_1(L_{0}(\xx)))}^*
                df_2(L_{0}(\xx))\rangle\,, \ \ \forall \xx\in
T^* \frakh^* \times \frakg^*_{red}; \label{eq:LL0}
\end{eqnarray}
\item Functions in $L_0^{*} I(L\frakg )$ provide a family
of Poisson commuting conserved quantities for $\HH_0$.
\item Under the Hamiltonian flow generated by $\HH_0$, we have
\begin{equation}
\frac{dL_{0}}{dt}=-[R_{0}^{*} (M_{0}), L_{0} ],
\end{equation}
where $M_{0}(q, p, s )(z)=L_{0}(q, p, s )(z)/z$, $\forall (q, p, s )
\in T^* \frakh^* \times  \frakg^*_{red}$.
\end{enumerate}
\end{thm}


\begin{rmk}
{\em 
\begin{enumerate}
\item   As in Remark  \ref{rmk:4.9}, we can derive a  St. Petersburg type
formula:
\begin{eqnarray}
\{L_0 (z) \  \stackrel{\ot}{,}\  L_0 (w)\}&=&-[\tilde{r}^{12}(q, z-w), \
 L^{1}(z)] +[\tilde{r}^{21}(q, w-z),  \ L^{2}(w)],
\end{eqnarray}
where $\tilde{r}(q, z)$ can be described  by an equation similar to 
Equation (\ref{eq:tr}).
 \item   Theorem \ref{thm:5.13} is presented  in a 
general Poisson setting. But it is clear that
 all the above claims are still  valid
when we  restrict the Hamiltonian $\HH_0$ as well as
other operators $L_{0}, \ R_{0}$ to a particular
symplectic leaf of  $T^* \frakh^* \times \frakg^*_{red}$.
\end{enumerate}
}
\end{rmk}

As in Remark \ref{rmk:5.1}, 
for any  $s \in \frakg^*_{red}$, we    write $s =\sum_{\alpha \in
\Delta} s_{\alpha} e^*_{\alpha}$ with $s_{\alpha_{i}}=1, \ 
i=1, \cdots , N$.
 Explicitly, the Hamiltonian  $\HH_0$ and the
Lax operators $L_0$ of  the integrable spin \cm systems
  on  $T^* \frakh^* \times  \frakg^*_{red}$
are given as follows:\\\\

{\bf I. {{ Rational case}}}

\be
&&\HH_{0} (q, p, s )=\half \sum_{i=1}^{N}p_{i}^{2}-\half \sum_{\alpha \in \Delta'}
\frac{s_{\al} s_{-\al} }{(\alpha , q)^2},\\
&&L_{0}(q, p, s )(z)=p+\frac{1}{z} \sum_{\alpha \in \Delta}
s_{\al} e_\al +\sum_{\alpha \in \Delta'}
\frac{ s_{\al}}{(\alpha , q)}e_\al ,
\ee
where  $\Delta' \subset \Delta$ is a set of roots
 closed with respect to the addition and multiplication by $-1$.\\\\

{\bf II. {{ Trigonometric  case}}}

\be
&&\HH_{0} (q, p, s)=  \half \sum_{i=1}^{N}p_{i}^{2} 
-\half \sum_{\al \in \Delta (\Pi'  )}(\frac{1}{\sin^{2}(\al, q)}\\
&&\ \ \ \ \ \ \ \ \ \ \ \ \ \ -\frac{1}{3})s_{\alpha} s_{-\al} 
 -\frac{5}{6}\sum _{\al \in \Delta-\Delta (\Pi'  )}s_{\alpha} s_{-\al}\\
&&L_{0}(q, p, s)(z)=p+  \sum _{\al \in \Delta (\Pi'  )}{ \mbox{sin}\, ((\al, q) + z)
\over
 \mbox{sin}\, (\al, q) \,  \mbox{sin}\, z}e^{\frac{1}{3}z (\al , q)}s_{\alpha} e_\al \\
&&\ \ \ \ \ \ \ \ \ \ + \sum _{\al \in \Delta_{+}-\Delta (\Pi'  )}\,{ e^{-iz}
\over
 \sin{ \,z}} e^{\frac{1}{3}z (\al , q)}s_{\alpha} e_\al \\
&&\ \ \ \ \ \ \ \ \ \  +\sum _{\al \in \Delta_{-}-\Delta (\Pi'  )}\,{ e^{iz}
\over
 \mbox{sin}\, z}e^{\frac{1}{3}z (\al , q)}s_{\alpha} e_\al .
\ee
Here $\Delta =\Delta_{+}\cup \Delta_{-}$
is a polarization of $\Delta$,
 $\Pi' $ is a  subset of the set of simple  roots,
and  $\Delta (\Pi'  )$
denotes the set of all roots which are linear
combinations of roots from $\Pi' $.\\\\

{\bf III. {{   Elliptic case }}}

\be
&&\HH_{0} (q, p, s)=  \half \sum_{i=1}^{N}p_{i}^{2}-\half
\sum_{\al \in \Delta }\calp ((\al , q )) s_{\alpha} s_{-\alpha}\\
&& L_{0} (q, p, s)(z)=p -\sum _{\al \in \Delta } l((\al , q ), z)
s_{\alpha} e_{\al}
\ee
where $\zeta (z)=\frac{\sigma' (z)}{\sigma (z)}$, $\calp (z)
=-\zeta' (z)$,
 $l(w, z)=-\frac{\sigma (w+z)}{\sigma (w)\sigma (z)}$, and
$\sigma (z)$ is the Weierstrass $\sigma$ function of periods
$2\omega_{1}, 2\omega_2$.\\

\begin{rmk}
{\em Let $\phi: \frakg \lon {\frak gl}(n, \complex )$ be a
representation of $\frakg$. Then it induces a representation
of $L\frakg$, which we denote  also by  the same
symbol. Let $A(q, p, \xi )=(\phi \smalcirc L_{0})(q, p, \xi )$,
where $L_{0}$ is the Lax operator of one of the 
integrable spin systems listed above. The we
have the spectral curve ${\cal C}: \ \mbox{det}(A(q, p, \xi )(z)-w)=0$,
which is preserved by the flow generated by the Hamiltonian 
$\HH_0$. The integrability of $\HH_0$ in the Liouville sense
on the symplectic leaves $T^* \frakh^* \times \calo_{red}$ of 
$T^* \frakh^* \times \frakg^*$ (of various dimensions) will be 
investigated in subsequent work.}
\end{rmk}
{\bf Example 5.1}
For $\frakg =\frak{sl}(3, \complex )$,   $\frakg^*_{red}$ can be
identified with  the affine subspace  consisting of matrices
of the form
$$s=\begin{pmatrix}
0&1&s_{13}\\
s_{21}&0&1\\
s_{31}&s_{32}&0
\end{pmatrix}.$$
The Poisson structure on $\frakg^*_{red}$ is given by
\be
\{s_{13}, s_{21}\}&=&1-s_{13}^{2}s_{21}\\
\{s_{13}, s_{31}\}&=& s_{13}(s_{21}-s_{32})\\
\{s_{13}, s_{32}\}&=& -1+s_{13}^{2}s_{32}\\
\{s_{21}, s_{31}\}&=& s_{21}(s_{32}-s_{13}s_{31})\\
\{s_{21}, s_{32}\}&=& s_{31}-s_{13}s_{21}s_{32}\\
\{s_{31}, s_{32}\}&=& s_{32}(s_{21}-s_{13}s_{31}).
\ee
If we consider the rational case with $\Delta' =\Delta$, then
$\HH_0$ and  $L_0$ are given as follows:
\be
&&\HH_{0} (q, p, s )=\half \sum_{i=1}^{2}p_{i}^{2}\\
&&\ \ \ \ \ \ \ \ \ \  -[\frac{s_{21}}{(q_{1}-q_{2})^2}+
\frac{s_{13}s_{31}}{(q_{1}-q_{3})^2}+\frac{s_{32}}{(q_{2}-q_{3})^2}]\\
&&L_{0} (q, p, s )=p+\frac{1}{z}s+\sum_{i\neq j}\frac{s_{ij}}{q_{i}-q_{j}}e_{ij},
\ee
where $s_{12}=s_{23}=1$ and $e_{ij}$ is the $3\times 3$ matrix
with a $1$ in  the $(i, j)$-entry and zeros elsewhere.\\\\



As a special case,  consider 
$\frakg={\frak{sl}}(N, \complex )$ and  identify $\frakg^*$ with
$\frakg$ using  the standard Killing form.  Let $\calo$ be the
adjoint  orbit through the point $\xi_{0} \in {\frak{sl}}({N}, \complex )$,
where $\xi_0$ is the  off-diagonal 
matrix  with all off-diagonal entries  equal to $m$ ($\neq 0$).
 It is simple to see that $\calo$ has 
dimension  $2(N-1)$, i.e., twice of the rank of the Lie algebra.
Hence, $\calo_{red}$ is  just one point. Consequently, if we
restrict the integrable spin systems  to this particular
symplectic leaf $T^*\frakh^* \times\{ {\mbox{pt}}\}$, we 
obtain the usual \CM systems with
 coupling constants  $m^2$.  Thus we
have recovered the following:

\begin{cor}
\cite{BAB1} \cite{BAB2}
The usual \CM  (rational, trigonometric and elliptic) systems 
 associated to the Lie
algebra ${\frak{sl}}({N}, \complex )$ admit a Lax operator
 $L_{0}: T^*\complex^{N-1}\lon L{\frak{sl}}({N}, \complex )$ and
an $r$-matrix formalism.
\end{cor}

\begin{rmk}
\label{rmk:5.3}
{\em  Note that the  above adjoint orbit $\calo$ is a semi-simple
orbit of ${\frak{sl}} ( N, \complex )$.  For other types of 
simple Lie algebras, unfortunately, there does not exist 
any semi-simple orbit of dimension equals to  twice the rank
of the Lie algebra \cite{J1} \cite{J2}.
On the other hand, there do exist minimal nilpotent orbits of
the correct dimension for ${\frak sp}(2N, \complex)$ \cite{J1}.
 However, the corresponding
Hamiltonian reduces to that of a free system  (without potential) in this
case.  In other words,  the integrable spin systems obtained
above do not contain  the usual \CM systems as subsystems
 for other types of simple Lie  algebras. }
\end{rmk}

\end{document}